\documentclass[12pt]{article}
\usepackage{amsmath,amsthm,verbatim,amssymb,amsfonts,amscd,diagrams, graphics}
\theoremstyle{plain}
\newtheorem{theorem}{Theorem}[section]
\newtheorem{corollary}[theorem]{Corollary}
\newtheorem{lemma}[theorem]{Lemma}
\newtheorem{proposition}[theorem]{Proposition}

\theoremstyle{definition}

\theoremstyle{remark}
\newtheorem{remark}[theorem]{Remark}

\newarrow{ul}---->
\newarrow{Backwards}<----

\newcommand{\td}[1]{\tilde{#1}}
\newcommand{\into}{\hookrightarrow}
\newcommand{\Z}{\mathbb{Z}}
\newcommand{\Q}{\mathbb{Q}}
\newcommand{\R}{\mathbb{R}}
\newcommand{\G}{\mathbb{G}}
\newcommand{\N}{\mathbb{N}}
\newcommand{\bd}{\partial}

\renewcommand{\H}{\mathbb H}

\newcommand{\mc}[1]{\mathcal{#1}}

\newcommand{\mf}{\mathfrak}

\newcommand{\im}{\text{im}}
\newcommand{\cok}{\text{cok}}

\renewcommand{\i}{\mathfrak i}
\renewcommand{\j}{\mathfrak j}

\begin{document}

\title{Singular chain intersection homology for traditional and super-perversities}
\author{Greg Friedman\\Yale University}
\date{July 12, 2004}
\maketitle
typeset=\today

\begin{abstract}
We introduce a singular chain intersection homology theory which generalizes that of King and which agrees with the Deligne sheaf intersection homology of Goresky and MacPherson on any topological stratified pseudomanifold, compact or not, with constant or local coefficients, and with traditional perversities or superperversities (those satisfying $\bar p(2)>0$). For the case $\bar p(2)=1$, these latter perversitie were introduced by Cappell and Shaneson and play a key role in their superduality theorem for embeddings. We further describe the sheafification of this singular chain complex and its adaptability to broader classes of stratified spaces.
\end{abstract}



\tableofcontents

\section{Introduction}

We fulfill two primary goals, each with the aim of providing some geometrical underpinnings of \emph{intersection homology theory}, which is an important tool in the study of stratified spaces with broad-reaching applications in algebraic geometry and representation theory (see \cite{Klei}):

The first goal is to generalize the singular chain intersection homology of King \cite{Ki} to provide a singular chain theory that yields the same intersection homology modules as the Goresky-MacPherson Deligne sheaf theoretic approach when the  perversity parameters are  \emph{superperversities}. Sheaf-theoretic intersection homology with this type of perversity parameter (satisfying $\bar p(2)\geq 1$) plays a crucial role, for example, in the superduality theorem of Cappell and Shaneson \cite{CS}. While the King singular chain approach can be used to define superperverse intersection homology   \emph{with  constant coefficients}, the resulting modules 
 do not agree with those obtained via the standard Goresky-MacPherson sheaf model of intersection homology, and, furthermore, the pre-existing singular chain theory does not readily extend to the important case of superperverse intersection homology with local coefficients. 
We here provide a version of the singular chain theory whose intersection homology modules  agree with those of the  sheaf model for both constant and local coefficients. 

Our second goal, which  we achieve simultaneously, is to demonstrate that our singular chain model also provides the correct sheaf theoretic intersectional homology modules for traditional perversities, even on \emph{non-compact topological pseudomanifolds}. It had been conjectured by King \cite{Ki} that such a  singular chain model should exist; we here provide the details.

The existence of a satisfactory singular chain intersection homology theory  enables the extension of superperverse intersection homology to more general filtered spaces, such as Quinn's manifold weakly stratified spaces \cite{Q1,Q2}), and, in fact, provides a reasonable (though non-axiomatic) way to \emph{define} intersection chain sheaves on such spaces. See Remark \ref{R: extension}, below, for more details.

\subsection{Background}\label{S: background}

 In order to clarify the goals of this paper, some historical remarks are in order to place our results in context.

\emph{Intersection homology}, as first introduced by Goresky and MacPherson \cite{GM1}, was  defined initially only on compact \emph{piecewise linear (PL) stratified pseudomanifolds} (a detailed definition of these and other relevant spaces is provided in Section \ref{S: spaces}, below). Working with fixed stratifications, intersection homology was defined as the homology of an  \emph{intersection chain complex}, given in terms of simplicial chains satisfying certain allowability conditions, restricting the dimensions of their intersections with the various strata. These allowability restrictions were determined by a fixed traditional \emph{perversity} function, i.e. a function $\bar p$ from $\N=\{0,1, 2,\ldots\}$ to $\N$ such that $\bar p(0)=\bar p(1)=\bar p(2)=0$ and $\bar p(k)\leq \bar p(k+1)\leq \bar p(k)+1$ for each $k\geq 2$. For technical reasons, these simplicial chains actually live in the direct limit over all compatible triangulations of the pseudomanifold, although it was later shown by Goresky and MacPherson in the appendix to \cite{MV86} that  the definition could be made with respect to a fixed triangulation provided it is flaglike with respect to the stratification (in particular, choosing any triangulation compatible with the stratification and then subdividing barycentrically provides a sufficient triangulation). In this incarnation, all chains have compact supports. This context was sufficient for Goresky and MacPherson to achieve their famed duality result that if $X$ is a compact $n$-dimensional stratified pseudomanifold and $\bar p$ and $\bar q$ are \emph{complementary perversities} satisfying $\bar p(k)+\bar q(k)=k-2$, then there is an intersection pairing $I^{\bar p}H_i(X)\otimes I^{\bar q}H_{n-i}(X)\to \Z$ that is nondegenerate upon tensoring all groups with the rationals $\Q$. 

In order to show that the intersection homology groups they had constructed are topological invariants and thus, in particular, independent of choice of stratification, Goresky and MacPherson then turned to sheaf theory in \cite{GM2}. Here it is shown that, on a PL stratified  pseudomanifold (not necessarily compact), if one instead begins with \emph{locally-finite} simplicial intersection chains, one can form a complex of sheaves of germs of intersection chains. This sheaf is based upon the presheaf whose sections over the open subspace $U$ are the locally-finite intersection chains on $U$. The intersection homology groups now appear as the hypercohomology of this complex of sheaves. On a compact pseudomanifold one obtains exactly the same groups as before, while on a non-compact pseudomanifold one attains a Borel-Moore type version of intersection homology theory (though the compact theory can also be recovered via hypercohomology with compact supports). In addition, this sheaf complex  is quasi-isomorphic to one suggested by Deligne whose quasi-isomorphism type can be completely described axiomatically and, even better, by a set of axioms independent of the particular choice of stratification. These properties were then utilized to demonstrate topological invariance of the intersection homology groups, which can now be recognized as the hypercohomology of this Deligne sheaf. 

The sheaf-theoretic formulation of intersection homology theory allowed the introduction of two important generalizations:
\begin{enumerate}
\item The Deligne sheaf can be defined  on \emph{topological} stratified pseudomanifolds and the axiomatic characterizations continue to demonstrate topological invariance of intersection homology in this setting

\item The Deligne sheaf can also be constructed beginning with a system of local coefficients which needs only be defined on the top dense manifold stratum of the pseudomanifold. This local coefficient theory gives up some of the topological invariance  (though we maintain such invariance with respect to stratifications compatible with the domain of definition of the coefficient system - see \cite{Bo}). Nonetheless, it provides a much richer theory allowing one, for example, to obtain invariants of embeddings of pseudomanifolds (see \cite{CS,GBF2}).
\end{enumerate}
As is well known (see, e.g. \cite{Klei} for a historical survey), the sheafified version of intersection homology theory has gone on to become an important tool not only in topology, but in algebraic geometry, representation theory, and the general theory of self-dual and perverse sheaves. 

A proof of the topological invariance of intersection homology not involving sheaves was then given by King in \cite{Ki}. This was done by introducing a chain complex of compactly supported  \emph{singular intersection chains}, again defined by their allowability with respect to a perversity function. With this definition, intersection homology \emph{with fixed coefficients} can, in fact, be defined on any filtered space, and King demonstrated topological invariance on the class of \emph{locally conelike topological stratified sets}. He also proved this invariance without requiring that $\bar p(0)=\bar p(1)=\bar p(2)=0$, though the other condition on perversities must be maintained. King notes that on a compact PL pseudomanifold, this singular chain theory agrees with the original PL theory of Goresky and MacPherson \cite{GM1}, and hence for traditional perversities ($\bar p(0)=\bar p(1)=\bar p(2)=0$), it also agrees with the sheaf theory on such spaces. He then conjectured the possibility of modifying the singular chain theory to obtain a complex of sheaves satisfying the axioms of the Goresky-MacPherson-Deligne sheaf and thus demonstrating that the singular chain and sheaf approaches provide the same theory in the \emph{topological category} even on non-compact spaces. We provide such a sheaf complex, built from singular chains, in this paper.

At this point, it is necessary to say a few words about \emph{superperversities}, those perversities for which $\bar p(2)\geq 1$.
The case $\bar p(2)=1$, for example, plays a key role in the intersection homology superduality theorem of Cappell and Shaneson \cite{CS}, which generalize Milnor's theorem for duality of infinite cyclic covers \cite{M68}.
For such perversities \emph{and for constant coefficients}, there is no difficulty in extending the definition of either the Goresky-MacPherson-Deligne sheaf or the King singular intersection chain complex. However, the resulting theories no longer agree, even on compact stratified PL pseudomanifolds. This can be observed as follows: In \cite{CS}, using the sheaf theoretic version of the theory, Cappell and Shaneson demonstrate that if $K$ is a locally-flat knot in $S^n$, $C$ is the complement of an open tubular neighborhood of $K$, and $\Lambda$ is a certain local coefficient system defined on the complement of $K$, then for $\bar p(2)=1$, $I^{\bar p}H_i(S^n;\Lambda)\cong H_i(C,\bd C;\Lambda)$. In particular, $I^{\bar p}H_0(S^n;\Lambda)=0$. This is impossible in any geometrically defined non-relative compact chain theory since any single point $0$-simplex cannot bound. 
It is further demonstrated by the author in \cite{GBF2} that even on $S^2$ with constant $\Z$ coefficients and a stratification $*\subset S^2$, the PL intersection chain sheaf does not satisfy the Deligne sheaf axioms if we use superperversities.  Hence there is a need to reconcile the geometric theory with the sheaf-theoretic one. In the local coefficient case, there is also the need to define exactly what one means by a geometric intersection chain complex, since now one must contend with $i-1$ or $i$ faces of $i$ simplices that can lie in the singular locus in which no coefficient is defined.

In \cite{HS91}, Habegger and Saper, also working in the category of locally conelike topological stratified spaces, presented a generalization of the Deligne sheaf construction which on PL pseudomanifolds and with constant coefficients provides hypercohomology modules that agree with the PL chain definition of intersection homology even for superperversities. They generalized further to \emph{codimension $\geq c$ intersection cohomology theories} defined with coefficients in a sheaf constructible with respect to some stratification of the space. Thus one can say that Habegger and Saper provide a sheaf version of the PL chain theory. We will take somewhat the opposite tack in finding a singular chain theory that when ``sheafified'' provides the same results as the pre-existing sheaf theory on \emph{topological} pseudomanifolds. 

We can also extend to yet a broader class of spaces. In \cite{Q2}, Quinn extended the study of constant coefficient compact singular chain intersection homology to \emph{manifold weakly stratified spaces} and demonstrated independence of stratification. On such spaces, local properties are specified not by topological conditions but by homotopy data. Hence singular chains seem to present a more natural approach then sheaf theory, as compactly supported intersection homology is a stratified homotopy type invariant (a proof of this long-standing folk theorem is provided by the author in \cite{GBF3}). Thus our approach seems to be a reasonable candidate for extending superperverse intersection homology  to such spaces so as to obtain a theory most closely resembling the sheaf version of the theory on pseudomanifolds. In fact, our approach extends easily to any filtered space, though of course we make no claim to topological invariance on such general spaces.

To summarize then, the singular chain intersection homology theory presented in this paper achieves the following goals:

\begin{itemize}
\item We provide a singular chain intersection homology theory that is well-defined both for traditional perversities and for more general perversities (including superperversities) on compact or non-compact spaces and with constant coefficients or local coefficients. 

\item On paracompact topological pseudomanifolds (compact or not, constant or local coefficients) and for traditional perversities, our intersection homology modules agree with those obtained by the Goresky-MacPherson-Deligne sheaf process.

\item On paracompact topological pseudomanifolds (compact or not, constant or local coefficients) and for superperversities satisfying $\bar p(2)=1$, our intersection homology modules agree with those utilized by Cappell and Shaneson in their superduality theorem.

\item There is a reasonable geometric extensions of our singular chain theory to any filtered space, and furthermore, on paracompact filtered Hausdorff spaces, our theory ``sheafifies'' in the sense that the intersection homology modules can be described both in terms of homology modules of a singular chain complex or as the hypercohomology modules of a homotopically fine sheaf. Thus we provide a (non-axiomatic) sheaf model for intersection homology on such spaces. 
\end{itemize} 

\subsection{Outline}

In Section \ref{S: chains} - \textbf{Chains}, we introduce our version of the singular intersection chain complex and study its properties. Section \ref{S: spaces} contains the  definitions of the types of spaces we will consider, and Section \ref{S: sing chains}contains  some notation concerning singular chain complexes.  In Section \ref{S: IC}, we first recall the construction of the singular intersection chain complex of King \cite{Ki} and then discuss the generalizations necessary to work with local coefficients. 
In Subsection \ref{S: G0}, we introduce a stratified coefficient system $\mc G_0$ determined by a given system of local coefficient $\mc G$ on the complement of the singular locus of a filtered space, and we define intersection homology with coefficients in $\mc G_0$.  It will be intersection homology with coefficients in $\mc G_0$  that allows us to obtain intersection homology modules isomorphic to superperverse Deligne-sheaf intersection homology modules with coefficients in $\mc G$. Section \ref{S: IC} also contains subsections on relative intersection homology (Section \ref{S: relative}), the stratified homotopy type independence of compactly supported intersection homology (Section \ref{S: sphe}), the stratification dependence of superperverse intersection homology (Subsection \ref{S: dep}), and intersection homology with more general ``loose'' perversities (Section \ref{S: codim 1}). 

In Section \ref{S: subdivision}, we show that the intersection homology class of an allowable chain is invariant under appropriate subdivisions, while in Section \ref{S: excision}, we demonstrate excision for intersection homology by proving the key result (Proposition \ref{P: U small}) that the singular intersection chain complex is chain homotopy equivalent to the intersection complex with supports in a locally-finite cover. This proposition also plays an important role in the proof, in Section \ref{S: sheaves}, that the sheaf built from our singular intersection chains is homotopically fine. Finally, Section \ref{S: compute} contains computations of the intersection homology of products, cones, and distinguished neighborhoods.

In Section \ref{S: sheaves} - \textbf{Sheaves}, we turn to sheaf theory and show that our singular intersection chain complex sheafifies to a sheaf whose hypercohomology agrees with both the homology of the singular chain complex and, on topological pseudomanifolds, the hypercohomology of the Goresky-MacPherson-Deligne sheaf. Section \ref{S: sheafify} contains the basic construction and  a demonstration that the sheaf we obtain is homotopically fine. In Section \ref{S: restriction}, we study the behavior of our intersection chain sheaf under restrictions to subspaces; in particular, we show that the restriction to an open subspace is quasi-isomorphic to the intersection chain sheaf of that subspace. Section \ref{S: axioms} contains the verification  of the agreement of the singular chain intersection  homology  the Deligne sheaf hypercohomology on pseudomanifolds. Lastly, in Section \ref{S: PL}, we indicate how superperverse intersection homology may be computed on PL pseudomanifolds via direct use of simplicial chains.

\section{Chains}\label{S: chains}

\subsection{Filtered spaces and stratified topological pseudomanifolds}\label{S: spaces}

In this section, we recall the definitions of the spaces in which we will be most interested.

A \emph{filtered space} is a topological space, $X$, together with a collection of
closed subspaces
\begin{equation*}
\emptyset = X^{-1}\subset X^0\subset X^1 \subset \cdots \subset X^{n-1}\subset X^n=X.
\end{equation*}
If we want to emphasize both the space and the filtration, we will refer to the
filtered space $(X,\{X^i\})$. Note that $X^i=X^{i+1}$ is possible. We will
refer to $n$ as the (filtered) \emph{dimension} of $X$ and to $X^{n-k}$ as the \emph{$n-k$
skeleton} or the \emph{codimension $k$ skeleton}. The sets $X_i=X^i-X^{i-1}$ are the \emph{strata} of $X$. We call a space either \emph{unfiltered} or \emph{unstratified} if we do not wish to consider any filtration on it (equivalently, $X=X^n$ and $X^i=\emptyset$, $i<n$).  

If $X$ is a filtered space, there is a canonical filtration of $X\times \R^k$ by $(X\times \R^k)^{i+k}=X^i\times \R^k$. If $cX$ denotes the open cone $X\times [0,1)/(x,0)\sim(y,0)$, there is a canonical filtration of  $cX$ so that $(cX)^{i+1}=X^i\times (0,1)$ for $i\geq 0$ and $(cX)^0$ is the cone point. If $X=\emptyset$, then by definition $c(X)$ is a point $x$ stratified as $c(X)=(c(X))^0=x$. Unless otherwise specified, all cones and products with $\R^k$ of filtered spaces are assumed to be given these canonical filtrations. Note that this dimensional indexing differs from the more codimensional indexing of \cite{GBF3}.

A filtered Hausdorff space $X$ is an \emph{$n$-dimensional topological stratified  pseudomanifold} if 
\begin{enumerate}
\item $X^{n-1}=X^{n-2}$, in which case $X^{n-2}$ is referred to as the \emph{singular locus} and denoted by $\Sigma$,

\item $X_k=X^{k}-X^{k-1}$ is either a topological manifold of dimension $k$ or it is empty,

\item $X-X^{n-2}=X-\Sigma$ is dense in $X$,

\item for each point $x\in X_{n-k}=X^{n-k}-X^{n-k-1}$, there exists a \emph{distinguished neighborhood} $N$ of $x$ such there is a compact topological stratified pseudomanifold $L$ (called the \emph{link} of the component of the stratum $X_{n-k}$), a filtration 
\begin{equation*}
L=L^{k-1}\supset  \cdots \supset L^0\supset L^{-1}=\emptyset,
\end{equation*}
and a homeomorphism
\begin{equation*}
\phi: \R^{n-k}\times c(L)\to N
\end{equation*}
that takes $\R^{n-k}\times c(L^{k-j-1})$ onto $X^{n-j}$. 
\end{enumerate}

A filtered Hausdorff space $X$ is an \emph{$n$-dimensional piecewise-linear (PL) stratified  pseudomanifold} if $X$ is a PL space, each $X^i$ is PL subspace, and, in the preceding definition, we replace topological manifolds with PL manifolds and homeomorphisms with PL homeomorphisms. 

A space is called simply an \emph{$n$-dimensional topological pseudomanifold} or an \emph{$n$-dimensional PL pseudomanifold} if it can be endowed with the structure of a, respectively, topological or PL stratified pseudomanifold. 
Intersection homology is known to be a topological invariant of such spaces; in particular, it is invariant under choice of  stratification (see \cite{GM2, Bo, Ki}).

If $X$ and $Y$ are two filtered spaces, we call a map $f:X\to Y$
\emph{stratum-preserving} if the image of each component of a stratum of $X$
lies in a stratum of $Y$ (compare \cite{Q1}). In general, it is not required
that strata of $X$ map to strata of $Y$ of the same (co)dimension. 
However, if $f$ preserves codimension then $f$ will induce a well-defined map on intersection homology (see \cite[Prop. 2.1]{GBF3}).  
We call $f$
a \emph{stratum-preserving homotopy equivalence} if there is a stratum-preserving
map $g:Y\to X$ such that $fg$ and $gf$ are stratum-preserving homotopic to the
identity (where the filtration of $X\times I$ is given by the collection
$(X\times I)^i=X^{i-1}\times I$). We will sometimes denote the stratum-preserving homotopy
equivalence of $X$ and $Y$ by $X\sim_{sphe}Y$ and say that $X$ and $Y$ are
\emph{stratum-preserving homotopy equivalent} or \emph{s.p.h.e.}
Stratum-preserving homotopy equivalences induce intersection homology isomorphisms \cite{GBF3}.

\subsection{Singular chains and related concepts}\label{S: sing chains}

Our most basic terms concerning singular chains should correspond to the standard concepts (see, e.g., \cite[\S29]{MK}):

Let $\Delta^i$ denote the standard affine $i$-simplex. A \emph{singular $i$-simplex} in a space $X$  is a continuous map $\sigma: \Delta^i\to X$. The image of a singular simplex is also referred to as its \emph{support} $|\sigma|$. \emph{Coefficients} of singular simplices will be defined in more detail below in Section \ref{S: IC}, dependent upon the particular homology theory under consideration, but all coefficients will  either be elements of a fixed module or certain lifts to bundles of modules. A \emph{finite $i$-chain} $\xi$ is a finite linear combination $\sum_j n_j\sigma_j$ of singular $i$-simplices $\sigma_j$ together with their coefficients $n_j$. A \emph{locally-finite $i$-chain} is a perhaps infinite formal sum $\sum n_j\sigma_j$ with the restriction that every point in $X$ possesses a neighborhood such that all but a finite number of simplices with support intersecting the neighborhood have $0$ coefficient (or equivalently are omitted from the sum). The \emph{support} of a chain $|\xi|$ is the union of the supports of its simplices with non-zero coefficients; the support of a finite chain is always compact. 

The collections of finite or infinite $i$ chains form abelian groups in the usual way. With coefficient system $\mc{G}$ these will be denoted $C^{c}_i(X;\mc{G})$ or  $C^{\infty}_i(X;\mc{G})$, respectively. For statements that hold for both finite and locally-finite chain groups, we will generally use the generic notation $C_i(X;\mc {G})$.   Boundary homomorphisms are given by the usual formula. 

If $\xi=\sum n_j\sigma_j$ is a chain and $\tau$ is a face of a singular simplex $\sigma_j$ such that $n_j\neq 0$  then by the \emph{star} of $\tau$ we mean the  chain $\sum n_k\sigma_k$, where the sum is taken over all $\sigma_k$  that have $\tau$  as a face.

\subsection{Intersection homology}\label{S: IC}

Here we review the basic definitions of singular intersection homology and introduce the coefficient system $\mc{G}_0$.

A \emph{ perversity} is a function $\bar p:\Z^{\geq 1} \to \Z$ such that $\bar p(k)\leq \bar p(k+1)\leq \bar p(k)+1$. A perversity is \emph{traditional} if $\bar p(1)=\bar p(2)=0$; these are the perversities originally employed by Goresky and MacPherson in their first definition of intersection homology \cite{GM1}. A perversity  is a \emph{superperversity} if  $\bar p(2)>0$. Study  of  nontraditional perversities occurs in a variety of sources, e.g. \cite{Ki, HS91}.   Our specific interest in superperversities stems from the key role of superperversities with $\bar p(2)=1$ in the superduality theorem of Cappell and Shaneson \cite{CS}. \emph{Unless otherwise specified, all perversities will be either traditional perversities or superperversities.} See Section \ref{S: codim 1}, below, for discussion of more general perversities. 

Singular intersection chains on a filtered space $X$ were first studied by King in \cite{Ki}. For a constant coefficient module $G$, the intersection chain complex $I^{\bar p}C^c_*(X;G)$ is a subcomplex of $C^c_*(X;G)$ defined as follows:  An $i$-simplex is called \emph{$\bar p$ allowable} if $\sigma^{-1}(X_{n-k})$ is contained in the $i-k+\bar p(k)$ skeleton of $\Delta^i$, and an $i$-chain $\xi$ is $\bar p$ allowable if each singular simplex in $\xi$ and $\bd \xi$ is $\bar p$ allowable. N.B. any singular simplex with $0$ coefficient is not considered in deciding allowability of a chain.  The intersection homology groups of King are then defined by $I^{\bar p}H^c_*(X;G)=H_*(I^{\bar p}C^c_*(X;G))$. 

An equivalent formulation of coefficients in the constant  module $G$ is to consider each coefficient $n_i$ of a simplex $\sigma_i$ as a lift of $\sigma_i$ to the $n_i$ section of the trivial bundle of modules $X\times G$ over $X$. Of course this approach readily generalizes to any bundle of coefficient modules $\mc G$ defined over $X$ - a coefficient of $\sigma$ is a lift of $\sigma$ to the bundle $\mc G$ and group operations are carried out continuously stalk-wise. This is the approach to homology with local coefficients espoused, for example, in Hatcher \cite{Ha}. 

Now suppose, however, that $X$ is a topological pseudomanifold and that $\bar p$ is a traditional perversity.
It is well known \cite{GM2} that to define a $\bar p$ intersection homology theory on $X$ with local coefficients, it is only necessary to specify a coefficient system $\mc G$ on $X-\Sigma$, the complement of the singular locus.
From the simplicial or singular chain point of view, this is essentially due to the fact that the allowability conditions on simplices in intersection chains prevent  them from intersecting the singular locus except in their codimension-two skeleta, so  boundary maps remain well-defined on coefficients. In this case, intersection homology need no longer be a topological invariant, although it will be for restratifications that are properly adapted to the coefficient system (see \cite[\S V.4]{Bo} for details). 

Under the assumptions of the preceding paragraph, we redefined the coefficient of a simplex slightly in \cite{GBF3} to be a lift only of   $ \Delta^i-\Delta^{i,i-2}$, the complement of the $i-2$ skeleton of $\Delta^i$. This was sufficient to obtain a well-defined intersection chain complex since the allowability conditions force each $\sigma(\Delta^i)$ to intersect $\Sigma$ only in $\sigma(\Delta^{i,i-2})$. Hence it is possible to lift $ \Delta^i-\Delta^{i,i-2}$ to $\mc G$ and to take boundaries in a reasonable way: for each $i-1$-face of each singular $i$-simplex in a prospective chain, we can lift at least its interior to $\mc{G}$ and perform the standard group operations  in this interior.

However, as noted by the author in \cite{GBF2}, it is not clear for superperversities how to define singular intersection homology with local coefficients defined only on $X-\Sigma$. If $\bar p(2)=1$,  the codimension 1 faces of simplices  may now dip into the singular locus $\Sigma$, and for higher superperversities, entire simplices may be allowed in $\Sigma$. 
In order to remedy this situation, it is necessary to extend the coefficients in some way into the singular set along with the chains. For this we utilize a stratified coefficient system $\mc G_0$, defined below. It is singular chain intersection homology with coefficients in $\mc G_0$ that will eventually recover for us sheaf  intersection homology with coefficients in $\mc G$. (We note that in the sheaf version of the theory, there is no difficulty in  extending to non-traditional perversities - one just follows the Deligne process \cite{GM2,Bo, CS}.)

\subsubsection{The coefficients $\mc{G}_0$.}\label{S: G0} Suppose $X$ is a filtered space and that a local coefficient system $\mc G$ of $R$ modules is given on $X-X^{n-1}$. We define $\mc G_0$ to consist  of the pair of coefficient systems defined by $\mc G$ on $X-X^{n-1}$ and the constant $0$ system on $X^{n-1}$. Given a singular simplex $\sigma: \Delta \to X$, a coefficient lift of $\sigma$ is then defined by a lift of $\sigma|_{ \sigma^{-1}(X-X^{n-1})}$ to $\mc G$ together with the trivial ``lift'' of $\sigma|_{ \sigma^{-1}(X^{n-1})}$ to the $0$ section $X^{n-1}\times 0$ of the $0$ system on $X^{n-1}$. A coefficient of a simplex $\sigma$ is considered to be the $0$ coefficient if it maps each point of $\Delta$ to the $0$ section of one of the coefficient systems. Note that if $\sigma^{-1}(X-X^{n-1})$ is path connected, then a coefficient lift of $\sigma$ to $\mc G_0$ is completely specified by the lift at a single point of $\sigma^{-1}(X-X^{n-1})$ by the lifting extension property for $\mc G$. In particular, for any traditional perversity or superperversity with $\bar p(1)=0$, the image of the interior of any allowable simplex $\sigma$ must lie in $X-X^{n-1}$ and thus $\sigma$ will have non-trivial coefficient lifts to $\mc G$ over, at least, this interior (in fact, the set of lifts of $\sigma$ is in bijective correspondence to the set of elements of the stalk $G$ of $\mc G$ over any point in the image of this interior).

Now we can define the intersection chain complex $I^{\bar p}C_*^c(C;\mc G_0)$. 
In this case  we still define a singular simplex $\sigma:\Delta^i\to X$ to be allowable if $\sigma^{-1}(X_{n-k}-X_{n-k-1})$ is contained in the $i-k+\bar p(k)$ skeleton of $\Delta^i$, but we define a coefficient $n$ of $\sigma$ in $\mc G_0$ to be a lift  of $\sigma^{-1}(X-X^{n-1})$ to $\mc G$, while the set $\sigma^{-1}(\Sigma)$ carries the zero coefficient lift over $X^{n-1}$. The boundary of $n\sigma$ is given by the usual formula $\bd(n\sigma)=\sum_j (-1)^j(n\circ \i_j)(\sigma\circ \i_j)$, where $\i_j:\Delta^{i-1}\to \Delta^i$ is the $j$th face map. Here $n\circ \i_j$ should be interpreted as the restriction of $n$ to the $j$th face of $\sigma$, restricting the lift to $\mc G$ where possible and restricting to $0$ otherwise. Of course the sign $(-1)^j$ multiplies the coefficient stalkwise. The boundary operator is extended from simplices to chains by linearity. $I^{\bar p}C_*^c(C;\mc G_0)$ is then the complex of finite chains $\xi=\sum n_i\sigma_i$ such that each simplex in $\xi$ and $\bd \xi$ is $\bar p$ allowable. Again, we note that any simplex with a zero coefficient is removed from the chain (or its boundary chain) and  is not taken into account for allowability considerations. In particular, a boundary face of a simplex in an allowable chain may not be allowable, so long as this boundary face ``cancels out'' of the boundary of the chain.

We should check that $I^{\bar p}C_*^c(X;\mc G_0)$ forms a legitimate chain complex, i.e. that the boundary of an allowable chain is allowable and that $\bd^2=0$. 
So let $\xi\in I^{\bar p}C_i^c(X;\mc G_0)$. Since $\xi$ is allowable, we know that the $i-1$ simplices in $\bd \xi$ with non-zero coefficients  must be allowable. 
So to show that $\bd \xi$ is allowable, it suffice to show that $\bd^2 \xi=0$, which is certainly allowable. But of course the operator $\bd^2$ is $0$ by purely formal considerations, stemming from its definition. 

So the allowable chains with respect to either a traditional perversity or a superperversity $\bar p$ and coefficient system $\mc G_0$ do form a chain complex $I^{\bar p}C^c_*(X;\mc G_0)$, and we denote its homology $I^{\bar p}H^c_*(X;\mc G_0)$. \emph{N.B. Throughout the paper, $I^{\bar p}H_*$ or $IH_*$ will always refer to the homology of the intersection chain complex defined in this section. Goresky-MacPherson intersection homology will be denoted in later sections as the hypercohomology of the Deligne sheaf $\H^*(\mc P^*)$.}

The above discussion applies equally well to locally-finite, perhaps infinite, chain complexes,  and so we can also define $I^{\bar p}C^{\infty}_*(X;\mc G_0)$ and $I^{\bar p}H^{\infty}_*(X;\mc G_0)$. We will continue to use the notation $I^{\bar p}C_*(X;\mc G_0)$ and $I^{\bar p}H_*(X;\mc G_0)$ in arguments that apply to both the complexes of finite and locally-finite chains. 

It is also important to notice that if $\bar p$ is a traditional perversity and $X^{n-1}=X^{n-2}$ (in particular if $X$ is a pseudomanifold), then for any allowable $i$-simplex $\sigma$, $\sigma^{-1}(\Sigma)$ will be contained in the codimension $2$ skeleton of $\Delta^i$, so, in particular, the interior of $\Delta^i$ and the interiors of its $i-1$ faces can always be given non-zero coefficient lifts. We then see that the chain complex obtained agrees completely with that for intersection homology with local coefficients in $\mc G$ as defined in \cite{GBF3}. Thus the following result follows immediately from the definitions, though we dub it a proposition simply to call attention to it.

\begin{proposition}
If $\bar p$ is a traditional perversity and $X^{n-1}=X^{n-2}$, then $I^{\bar p}C_*(X;\mc G_0)=I^{\bar p}C_*(X;\mc G)$. 
\end{proposition}

In particular, in all of our following computations involving $I^{\bar p}C_*(X;\mc G_0)$ and $I^{\bar p}H_*(X;\mc G_0)$, if $\bar p$ is a traditional perversity and $X^{n-1}=X^{n-2}$, we recover (or obtain!) results about $I^{\bar p}C_*(X;\mc G)$ and $I^{\bar p}H_*(X;\mc G)$. For this reason, and also because we will show that $I^{\bar p}H_*(X;\mc G_0)$ agrees with sheaf intersection homology on pseudomanifolds if $\bar p$ is a superperversity, we will primarily restrict discussion in this paper to coefficients $\mc G_0$, \emph{even if $\mc G$ can be extended  to higher codimension strata of $X$}. However, all results of this paper except for those involving explicit computations of intersection homology modules remain valid if we replace $\mc G_0$ by other appropriate extensions of $\mc G$, in particular if $\mc G=G$ is a constant coefficient module. We leave the necessary modifications of the computational results of Section \ref{S: compute} to the interested reader.

\begin{remark}
Our approach to superperverse intersection homology via a pair of coefficient systems is somewhat foreshadowed by that of MacPherson in his unpublished monograph \cite{MacP90}. His procedure is to define the intersection chain complex on Whitney stratified spaces via locally-finite chains on the complement of the singular locus $X-\Sigma$, which, in effect, gives $0$ boundaries to chains that ``fall off the end'' of $X-\Sigma$ into $\Sigma$. However, this approach is developed not with singular chains but with certain ``good'' classes of geometric chains; furthermore, even on compact spaces, there is the necessity of working with locally-finite infinite chains. Our intersection  chain complex  has the potential advantage of requiring the use only of finite singular chains on compact spaces, and our singular chains are actually allowed  to intersect the singular locus, if permitted to do so by the allowability conditions. Our approach also treats a broader class of spaces than considered in \cite{MacP90}.
\end{remark}

\begin{remark}
It is worth pointing out that, even with constant coefficients $G$, use of the coefficient system $G_0$, in which intersections of simplices with $\Sigma$ carry a formal $0$ coefficient, is \emph{not} the same as attempting to take relative intersection homology $I^{\bar p}H_*(X,\Sigma;G)$. For one thing, there is no such chain submodule as $I^{\bar p}C_*(\Sigma; G)$ as no allowable chains are contained entirely within $\Sigma$. This is also not the same as killing a homotopy equivalent neighborhood $N$ of $\Sigma$ by considering $I^{\bar p}H_*(X,N;\mc G)$. In this case, we could excise $\Sigma$ (see Lemma \ref{L: excision2} below), and the module would reduce to $H_*(X-\Sigma,N-\Sigma;  G)$. However, the computations of intersection Alexander polynomials of singular knots in \cite{GBF2} adequately demonstrate that these modules do not agree with the Deligne sheaf hypercohomology and hence, as will be seen, do not  agree with the intersection homology modules defined in this section.
\end{remark}

\subsubsection{Relative intersection homology.}\label{S: relative} If $U$ is a subspace of $X$, then we define $I^{\bar p}C_*(U_X;\mc G_0)$ to be the  subcomplex of $I^{\bar p}C_*(X;\mc G_0)$ consisting of allowable chains (in $X$!) with support in $U$. Note that for finite chains $I^{\bar p}C_*^c(U_X;\mc G_0)\cong I^{\bar p}C_*^{c}(U;\mc G_0)$, while the analogous statement does not hold for $IC^{\infty}_*$ since there may be locally-finite chains in $U$ that are not locally-finite chains in $X$ (e.g. they may accumulate at a point in $\bar U-U$). We define $I^{\bar p}C_*(X,U;\mc G_0)=I^{\bar p}C_*(X;\mc{G}_0)/I^{\bar p}C_*(U_X;\mc G_0)$. These chain complexes yield intersection homology modules $I^{\bar p}H_*(U_X;\mc G_0)$ and $I^{\bar p}H(X,U;\mc G_0)$. While these definitions hold formally for any subspace $U\subset X$, in applications involving stratified pseudomanifolds  (in which case the skeleton $X^k$ actually has dimension $k$), one often wants $U$ itself to be a stratified pseudomanifold. In these cases, one usually considers only open subsets $U$, in which case the restricted stratification does provide a pseudomanifold stratification.

\subsubsection{Stratified homotopy invariance.}\label{S: sphe} As for traditional intersection homology, $I^{\bar p}H_*(X;\mc G_0)$ is not a homotopy invariant of $X$. However, $I^{\bar p}H^c_*(X;\mc G_0)$ it is an invariant of stratum-preserving homotopy type. For traditional perversities, a proof is given in  \cite{GBF3} (though it was certainly a folk-theorem beforehand). This result easily carries over to the more general cases considered here (though only with compact supports).

\begin{lemma}\label{L: sphe}
$I^{\bar p}H^c_*(\quad;\mc G_0)$ is a stratum-preserving homotopy invariant, i.e. any stratum-preserving homotopy equivalence $f:X\to Y$ induces an isomorphism on intersection homology. More specifically, $I^{\bar p}H^c_*(Y;\mc G_0)\cong I^{\bar p}H^c_*(X;(f^*\mc G)_0)$.
\end{lemma}
\begin{proof}
The proof is essentially that presented in \cite{GBF3} for intersection homology with traditional perversity and a local coefficient system $\mc G$ over $X-\Sigma$. The modifications needed to handle the more general cases are minor. 
\end{proof}

\begin{corollary}
If $f: (X,A)\to (Y,B)$ is a stratum-preserving homotopy-equivalence of pairs, then $I^{\bar p}H^c_*(Y, B;\mc G_0)\cong I^{\bar p}H^c_*(X,A;(f^*\mc G)_0)$.
\end{corollary}
\begin{proof}
This follows from the preceding lemma and the five lemma applied to the induced map of long exact sequences. 
\end{proof}

\subsubsection{Stratification dependence of superperverse intersection homology}\label{S: dep}

It is crucial to note that, even for constant coefficients on pseudomanifolds, if $\bar p$ is a superperversity, then $I^{\bar p}H_*(X;\mc G_0)$  will not generally be a topological invariant, completely independent of the stratification of $X$. In fact, consider the sphere $S^n$, $n\geq 2$. If $\bar p$ is a traditional perversity, then $I^{\bar p}H_*(S^n;\Z_0)\cong I^{\bar p}H_*(S^n;\Z)\cong H_*(S_n; \Z)$ for any stratification of $S^n$. Similarly, if $\bar p$ is a superperversity and $S^n$ is given the trivial stratification with no strata of codimension greater than zero, $I^{\bar p}H_*(S^n;\Z_0)=I^{\bar p}H_*(S^n;\Z)\cong H_*(S_n; \Z)$. But now, let $\bar p$ be the superperversity $\bar p(k)=k-1$, and suppose that $S^n$ is filtered by $x\subset S^n$, for some point $x\in S^n$. Then $I^{\bar p}H_0(X;\mc G_0)=0$ since for any allowable $0$-simplex represented by the map $z\in S^n-x$ (thought of, of course, as a simplex $z:\Delta^0\to X-\Sigma$), there is an allowable $1$-simplex $\sigma\in C_1(S^n;\Z)$ such that $\bd \sigma=z-x\in C_0(X;\Z)$. However, as an element of $I^{\bar p}C_1(X;\Z_0)$, $\bd \sigma=z$.  So all intersection $0$-cycles bound, and $I^{\bar p}H_0(S^n;\Z_0)=0$. 

Similarly, by adding more points to the $0$-skeleton, we cause $I^{\bar p}H_1(S^n;\Z_0)$ to be a free abelian group of any rank, generated by $1$-cycles consisting of curves connecting points of the $0$-skeleton. 

This stratification dependence should not be a concern, however, as in most patterns of application, e.g. \cite{CS, GBF2, Ba02}, a certain stratification is either assumed or forced upon us, and we wish to use intersection homology to study the space together with its stratification. Furthermore, if $X$ is a topological pseudomanfiold, $I^{\bar p}H_*(X;\mc G_0)$ \emph{will} be independent of changes of stratification that fix the top skeleton $X^{n-1}$. In other words, the intersection homology modules $I^{\bar p}H_*(X;\mc G_0)$ will agree for any two stratifications which share the same $n-1$ skeleton $X^{n-1}$. This is proven by the author in \cite{GBF11} for superperverse sheaf theoretic intersection homology on topological pseudomanifolds, so it will follow for the singular theory once we show that the two theories agree in Section \ref{S: sheaves}, below (our proof here does not rely upon this stratification invariance).

\subsubsection{Loose perversities}\label{S: codim 1}

Throughout this paper, we principally limit ourselves to traditional or super-perversities, i.e. those for which $\bar p(1)=\bar p(2)=0$ or $\bar p(2)>0$, respectively. These are the cases of greatest historical interest. In this section, we briefly discuss intersection homology with ``looser'' perversities.

First, let us consider a fixed coefficient group $G$ on $X$. In this case, King \cite{Ki} first defined singular intersection chains on filtered spaces  for any \emph{loose perversity}. A loose perversity is \emph{any} sequence of integers $\bar p(1), \bar p(2), \ldots, \bar p(n))$. Note that we are free to ignore $\bar p(0)$ or simply to assume $\bar(0)=0$ since  setting $\bar p(0)>0$ has no added benefit (an $i$ simplex can't intersect $X-X^{n-1}$ in a $>i$ skeleton),  while $\bar p(0)<0$ leads to having no allowable chains that intersect $X-X^{n-1}$, in which case we could simply restrict to an intersection homology theory on $X^{n-1}$. One then defines allowable chains word-for-word as before but using a loose perversity parameter. This process can clearly be extended to include locally-finite intersection chains. These completely  general perversities  are rarely used in practice, however.

We could also treat intersection homology with loose perversities and local coefficients. Once again, if $\mc G$ is defined on $X-X^{n-1}$, we can extend it to $\mc G_0$ and define $I^{\bar p}C_*(X;\mc G_0)$ in the obvious way. However, more interesting situations can occur. For example, a loose perversity might make  it possible for an $i$-chain to intersect some lower strata in a significant way and others not at all. It would be interesting to study what happens if we define intersection chain complexes that leave the original coefficients $\mc G$ on some strata and add the $0$ coefficient system on others. One could also put different, but compatible, coefficient systems on each stratum. 

We will not treat such general theories in this paper, but we do note that either for constant coefficients or for coefficients $\mc G_0$ as defined previously, the results of this paper  hold for loose perversities, barring those involving specific computations of intersection homology modules in Sections \ref{S: compute} and \ref{S: axioms}.

\subsection{Subdivision}\label{S: subdivision}

In this section, we will show that intersection homology classes are preserved under suitably defined subdivisions of their representative chains. Of course this is well-known for, e.g., barycentric subdivisions of singular chains in ordinary homology (see \cite[\S 31]{MK}). We will require more general subdivisions, and we must verify that $\bar p$ allowability is preserved. We begin by considering what it should mean for a singular chain to have a subdivision. We proceed by defining \emph{singular subdivisions} of certain simplicial complexes, and then we use these model singular subdivisions to obtain subdivisions of singular chains.

If $\Delta^0$ is the positively-oriented simplicial $0$ simplex, then a \emph{singular subdivision} of $\Delta^0$ is just the singular $0$-simplex $\i: \Delta^0 \to \Delta^0$. We let $-\i$ be the singular subdivision of $\Delta^0$ considered with the opposite orientation (so $-\i$ is a singular subdivision of $-\Delta^0$). 

Next, let $\Delta^{i}$, $i>0$, be the standard model oriented $i$-simplex. Let $\Delta'$ be a simplicial subdivision of $\Delta^{i}$, and let $\{\Delta^{i}_j\}_{j\in J}$ be the collection of simplicial $i$-simplices in the subdivision, all oriented consistently with $\Delta^{i}$.
A \emph{singular subdivision} of $\Delta^i$ corresponding to $\Delta'$ is a singular chain $s=\sum_j (-1)^{\text{sgn}(\i_j)} \i_j$ where 

\begin{enumerate}
\item each $\i_j: \Delta^i \to \Delta^i_j$ is a linear injection,

\item \label{I: orient} $\text{sgn}(\i_j)$ is $0$ if $\i_j$ is orientation preserving and $1$ if it is orientation reversing,

\item \label{I: bd} $\bd s$ is a singular subdivision of $\bd \Delta^i$ corresponding to the  subdivision induced by $\Delta'$ (in particular the support $|\bd s|\subset \bd \Delta$ and both the singular and simplicial boundary terms should have compatible orientations with respect to the induced orientation in the same sense as in item \eqref{I: orient}).
\end{enumerate}

Note that $-s$ is a singular subdivision of $-\Delta^i$, meaning $\Delta^i $ with the opposite orientation but with the same simplicial subdivision.

Given a subdivision $\Delta'$ of $\Delta^i$, a singular subdivision can always be obtained as follows: first partially order the vertices of $\Delta'$ such that the ordering is a total ordering on any collection of vertices that span a simplex and such that it preserves the canonical ordering on the vertices of $\Delta^i$. Then to each simplex $\Delta_j=[v_0,\ldots, v_i]\in \Delta'$ with $v_0<\cdots<v_i$, assign the singular simplex $\i_j$ given by the simplicial map $\Delta^i=[0, \ldots, i] \to [v_0,\ldots, v_i]$ determined by $k\to v_k$. Then $s=\sum (-1)^{\text{sgn}(\i_j)}\i_j$. This assignment essentially gives the standard image of the orientation class $\Gamma$ for $\Delta'$ under the chain map  $\phi: \mf C_*(\Delta')\to \mf C_*'(\Delta^i)\to C_*(\Delta')$ from oriented to ordered to singular chains of $\Delta'$ that is used in the usual proof of equivalence of singular and simplicial homology (see \cite{MK}). The first two conditions of the definition are clearly satisfied, and since $\phi$ is a chain map, $\bd \phi(\Gamma)=\phi\bd (\Gamma)$, which implies condition \eqref{I: bd}.
Conversely, given a singular subdivision, it is easy to see that the subdivision determines such a partial  ordering of the vertices of $\Delta'$ - just order by the standard ordering on $\Delta^i$ under the homeomorphism $\i_j$. 

If $\sigma: \Delta^i \to X$ is a singular $i$-simplex, then a subdivision $\sigma'$ of $\sigma$ with respect to the subdivision $\Delta'$ of $\Delta^i$ is the singular chain $\sigma s$, i.e. if $s=\sum (-1)^{\text{sgn}(\i_j)} \i_j$, then $\sigma s=\sum (-1)^{\text{sgn}(\i_j)} \sigma\circ \i_j$. 
If $n$ is a coefficient of $\sigma$, we similarly define the singular chain with coefficients $n'\sigma'=\sum (-1)^{\text{sgn}(\i_j)} ( n\circ \i_j)\sigma\circ \i_j$.

Suppose now that $\xi$ is an $i$-chain of $X$, $\xi=\sum n_k\sigma_k$. For each $\sigma_k$, let $\Delta_k$  represent a copy of the standard model $i$-simplex so that $\sigma_k:\Delta^i_k\to X$. We say that subdivisions $\{\Delta_k'\}$ are compatible with respect to $\xi$ if  the following condition holds: suppose that $\sigma_k$ and $\sigma_l$ are singular simplices in $\xi$ with non-zero coefficients and that they have faces $\tau_k$ and $\tau_l$ such that $\tau_k=\tau_l$ as singular simplices (i.e. $\tau_k: \Delta^{i-1}\to \Delta^i \overset{\sigma_k}{\to} X$ equals $\tau_l: \Delta^{i-1}\to \Delta^i \overset{\sigma_l}{\to} X$). Then the induced subdivision $\tau_k'$ and $\tau_l'$ should agree as chains. Note that $k$ may equal $l$ so this condition may impose non-trivial relations among faces of the same $i$-simplex. Given such compatible subdivisions, we can form the chain $\xi'=\sum n'\sigma'$ and have $\bd \xi'=(\bd \xi)'$, where the latter term indicates the induced subdivision of $i-1$ chains in the boundary of $\xi$. We call $\xi'$ a subdivision of $\xi$. A subdivision of a finite (locally-finite) chain is itself finite (locally-finite). 

The standard example of  a subdivision $\xi'$ is given by the barycentric subdivision of singular chains (see \cite{MK}). In this case, there is a natural partial ordering on the vertices of the subdivided model simplices determined by the dimension of the face of which each vertex is a barycenter. The uniformity of the construction ensures compatibility among simplices in any chain. Similarly, we can find such natural orderings for generalized barycentric subdivisions, in which not every face is subdivided at each step. In this cases, it is only necessary to find a scheme by which compatibility among simplices is maintained. Such a procedure is used in the proof of Proposition \ref{P: hom fine} below.

It will often be convenient in what follows to identify the images of  the subdivision singular simplices under the linear injections $\i: \Delta^i\into \Delta^i$ with their corresponding model simplices. In other words, we sometimes identify the singular simplex $\sigma \i$ with $\sigma$ restricted to the image of $\i$, which will be some subsimplex $\delta\subset \Delta^i$. This often makes the wording more convenient in arguments where we must check allowability conditions. It should always be remembered though that the specification of a subdivided simplex requires not just a restriction of $\sigma$ but a precise specification of how the model simplex is identified with $\delta$.  

We also note for future use that the idea of a singular subdivision of a simplex $\sigma:\Delta^i\to X$ modeled upon some subdivision $\Delta'$ can be extended to define singular subdivisions of any dimensionally homogeneous polyhedral space  based upon some triangulation by oriented simplices. In particular, we will need below such singular triangulations of $\Delta\times [0,1]$, $\Delta\times [0,\infty)$,  $\Delta\times \R$, and $\Delta\times \R^k$. In each case, we begin with a simplicial triangulation of the space and then use some partial ordering on the vertices to determine a singular triangulation.

Of course it will be important to know that the subdivision of a $\bar p$ allowable chain remains allowable:

\begin{lemma}\label{L: sub a}
Let $\xi'$ be a subdivision of the $i$-chain $\xi\in I^{\bar p}C_i(X;\mc G_0)$. Then $\xi'\in I^{\bar p}C_i(X;\mc G_0)$.
\end{lemma}
\begin{proof}
Recall that the $\bar p$ allowability of $\xi$ means that each  $i$-simplex  $\sigma$ in $\xi$ with non-zero coefficient satisfies the property that  $\sigma^{-1}(X_{n-k}-X_{n-k-1})$ is contained in the $i-k+\bar p(k)$ skeleton of $\Delta^i$ and similarly each $i-1$ simplex in $\bd \xi$ satisfies the analogous property with $i-1-k+\bar p(k)$. Now $\xi'$ is composed of the singular  $i$ simplices of the form $\sigma\i_j$ where $\i_j:\Delta^i\to \Delta^i$ is linear and injective. We must determine if $(\sigma\i_j)^{-1}(X_{n-k}-X_{n-k-1})$ is contained in the $i-k+\bar p(k)$ skeleton of $\Delta^i$. But note that $\i_j^{-1}$ of the $r$ skeleton of $\Delta^i$ must lie in the $r$ skeleton of $\Delta^i$ since $\i_j(\Delta^i)$ is an $i$-simplex of a subdivision of $\Delta^i$. 
Thus since $(\sigma\i_j)^{-1}(X_{n-k}-X_{n-k-1})=\i_j^{-1}\sigma^{-1}(X_{n-k}-X_{n-k-1})\subset \i_j^{-1}(  \{i-k+\bar p(k) \text{ skeleton of }\Delta^i\})\subset \{i-k+\bar p(k) \text{ skeleton of }\Delta^i\}$, we see that each $\sigma\i_j$ is allowable. Since $\xi'$ is composed of $i$-simplices of this form, we see that all its $i$-simplices are allowable. 
Similarly, the simplices in $\bd \xi'$ are allowable since $\bd \xi'$ is a subdivision of $\bd \xi$, so the above arguments hold one dimension lower. 
\end{proof}

Using this lemma, we show that an intersection cycle and its subdivisions define the same intersection homology class.

\begin{proposition}\label{P: subdivision}
Let $\xi$ be a $\bar p$ allowable cycle representing an element of $I^{\bar p}C^{\infty}_i(X, U;\mc G_0)$, where $U$ is possibly empty. Then $\xi$ is intersection homologous to any subdivision $\xi'$, so  $\xi$ and $\xi'$ represent the same element of $I^{\bar p}H^{\infty}_i(X,U;\mc G_0)$. If $\xi$ is a finite chain, then the realizing homology can also be taken as finite, and, in particular, the same statements hold for $I^{\bar p}C^{c}_*(X,U;\mc G_0)$ and $I^{\bar p}H^{c}_*(X;\mc G_0)$.
\end{proposition}
\begin{proof}
We can construct the homology rather explicitly by constructing an allowable $i+1$ chain $D$ such that $\bd D=\xi'-\xi+E$, where $E$ is an allowable chain in $X$ with support in $U$. We follow a fairly standard prism construction.

Suppose that $\xi=\sum n_j\sigma_j$ and let $\Delta_j$ be the model simplex for $\sigma_j$. By definition of $\xi'$, $\xi'=\sum_j \sum_l (-1)^{\text{sgn}(\i_{jl})}( n_j\circ \i_{j,l})\sigma\circ \i_{j,l}$, where each $\i_{j,l}$ is a linear injection $\Delta^i\to \Delta^i_j$ determined by a simplicial subdivision $\Delta_j'$ of $\Delta^i_j$. 

We begin by triangulating the set $B=\amalg \Delta_j\times [0,1]$. Suppose that each $\Delta_j\times 0$ is triangulated as its own simplex and that each $\Delta_j\times 1$ is triangulated as per our given subdivision $\Delta'$. We want to extend this triangulation to the whole space. The simplest procedure is inductive on the dimensions of faces in $\amalg \Delta_j\times 0$: For each $0$ simplex $v$ in $\amalg \Delta_j\times 0$ add the corresponding $1$-complex $\bar c([v\times 0] \cup [v\times 1])$, where $\bar c$ represents the closed cone. Now for each $k$ simplex $w$ in $\amalg \Delta_j\times 0$, assume that we have a triangulation of $C=(w\times 0) \cup (w\times 1) \cup (\bd w\times [0,1])$ (this triangulation is determined by the standard triangulation of $\Delta\times 0$, the subdivided triangulation of $\Delta\times 1$, and by the induction hypothesis). Now triangulate $w\times [0,1]$ by taking the closed cone on $C$. Since the dimension of the chain  $\xi$ is finite, this process terminates with a triangulation of $B$. Notationally, we denote $B$ with this  triangulation as $\td B$ and its restriction to each $\Delta_j\times [0,1]$ by $\td \Delta_j$.  
We take the partial ordering on the vertices of $\td B$ as determined on each simplex by first following the partial orderings of each $\Delta_j'$ and $\Delta_j$  and then ordering the cone points $c_a$ by the stage of their addition.

We will next construct a singular chain $S_j$ that will serve as a generalized singular subdivision of the polyhedron $\Delta_j \times [0,1]$ adapted to the triangulation  $\td \Delta_j$. Each $i+1$ simplex in the triangulation of $\td \Delta_j$ can be written in the form $\pm[\tau^k, c_1, \ldots, c_{i+1-k}]$, where $\tau^k$ is a $k$-simplex of $\Delta_j \times \{0,1\}$ and the $c_a$ are the cone points of the construction. To each such $i+1$ simplex, we assign a singular simplex $\j$ determined by the partial ordering. Then we take $S_j=\sum (-1)^{\text{sgn}(\j)} \j$, where the sign is determined by whether or not $\j$ agrees with the orientation of $\Delta \times [0,1]$ induced by the standard orientations. Now take $S=\sum S_j$ as a singular subdivision of $\td B$. 
Equivalently,  let $\phi: \mf C(\td B)\to \mf C'(\td B)\to C(\td B)$ be the map  from simplicial to singular chains via ordered simplicial chains, and, as before, let $S=\phi(\Gamma)$, where $\Gamma$ is the fundamental chain for $\td B$ as determined from the orientations.  $S$ is precisely the chain we have described in detail. Since $\phi$ is a chain map,  the support of $\bd S$ is in $\bd (B\times [0,1])$, and, in fact,  $\bd S_j=s'-\text{id}_{\Delta_j} +E_j$, where $|E_j|\subset \bd \Delta_j\times [0,1]$.

Finally, let $p:B\times [0,1] \to B$ be projection and define $$D(\xi)=D(\sum n_j\sigma_j)=\sum_j \sum_l (-1)^{\text{sgn}(\i_{jl}} (n_jp\j_{jl})\sigma_jp\j_{jl},$$ where the sum in $l$ is over all $i+1$ simplices in the triangulation of $\Delta_j\times [0,1]$ and $\j_{jl}$ is the singular chain corresponding to the $l$th singular simplex. Using the above computation for $\bd S_j$ and the obvious compatibility of the subdivision, we see that $\bd D(\xi)=\xi'-\xi +E$, where $|E|\subset |\bd \xi|\subset U$. 

It is clear from this construction that if $\xi$ is finite or locally-finite then so is $D(\xi)$. It remains to check that $D(\xi)$ is allowable. 
Since  $\xi$ and $\xi'$ are allowable  it remains only to check that the $i+1$ simplices of $D(\xi)$ are allowable and that the $i$ simplices of $E$ are allowable. So let $\eta=\sigma p\j$ be a singular $i+1$ simplex of $D(\xi)$. We have   $\eta^{-1}(X_{n-k}-X_{n-k-1}) =\j^{-1}(\sigma^{-1}(X_{n-k}-X_{n-k-1})\times [0,1])$. Since $\sigma^{-1}(X_{n-k}-X_{n-k-1})\subset  \{i-k+\bar p(k) \text{ skeleton of }\Delta_j\}$, $\sigma^{-1}(X_{n-k}-X_{n-k-1})\times [0,1]$ must lie in the $(i-k+\bar p(k))+1$ skeleton of our subdivision of $\Delta_j\times [0,1]$, and hence it intersects only the $(i-k+\bar p(k))+1$ skeleton of $\j(\Delta)$. But this implies that $\eta$ is $\bar p$ allowable. $E$ is allowable by exactly the same arguments one dimension lower by using the allowability of $\bd \xi$ and $\bd \xi'$.
\end{proof}

\subsection{Excision}\label{S: excision}

Next, we need a proposition that shows that it is possible to break intersection chains into small pieces, at least up to chain homotopy. This mirrors the usual proof of excision for singular chains (see e.g. \cite{MK}) except that more care must be taken to ensure allowability at each step along the way.

Let $\mc{U}=\{U_k\}$ be a locally-finite open cover of $X$. We choose and fix a well-ordering on $\mc{U}$. Let $I^{\bar p}_{\mc{U}}C^c_*(X;\mc{G}_0)$ be the subcomplex of $I^{\bar p}C^c_*(X;\mc {G}_0)$ consisting of intersection chains $\xi$ that can be written as the finite sum of intersection chains $\xi=\sum \xi_j$ such that each $\xi_j$ has support in some  $U_k$. Let $\iota : I^{\bar p}_{\mc{U}}C^c_*(X;\mc G_0)\to I^{\bar p}C^c_*(X;\mc {G}_0)$ be the inclusion. We will see that this inclusion is a chain homotopy equivalence. 

\begin{remark}We cannot expect to obtain a similar statement concerning $I^{\bar p}C^{\infty}_*(X;\mc {G}_0)$ (at least if $|\mc U|=\infty$) since a chain composed of an infinite number of singular simplices cannot be written as a sum of an the arbitrary number of pieces. The sum that occurs in $\xi=\sum_jn_j\sigma_j$ is a formal one and does not correspond to group operations, which cannot be infinitely strung together. 
\end{remark}

\begin{proposition}\label{P: U small}
There exists a chain map $T: I^{\bar p}C^c_*(X;\mc {G}_0) \to I^{\bar p}_{\mc{U}}C^c_*(X;\mc {G}_0)$ and a chain homotopy $D$ from $\iota T$ to the identity. 
\end{proposition}

\begin{proof}
We define first a chain map $\hat T:C^c_*(X)\to {}_{\mc U}C^{c}_*(X)$ from singular chains with coefficients in $\Z$ to singular chains supports in $\mc U$. We will then use $\hat T$ to construct $T$. 

Throughout the proof, we fix a function  $\psi$ assigning to each singular simplex $\sigma: \Delta \to X$ with support contained in some element of $\mc U$ the smallest $U_k$ such that $|\sigma|\subset U_k$ (here we use the fixed ordering on $\mc{U}$). Note that $\psi$ depends only on $|\sigma|$, not on the specific map.

The definition of $\hat T$ proceeds by induction. Since $C^c_0(X)={}_{\mc U}C_0^c(X)$, we let $\hat T$ be the identity map on $C^c_0(X)$. 

Now suppose by induction hypothesis that we have defined $\hat T$ on $C^c_j(X)$ for all $j\leq i-1$ and that for each $j$-simplex $\tau$, the following conditions are satisfied:

\begin{enumerate}
\item $\hat T\tau$ is a subdivision of $\tau$.

\item $\hat T$ is a chain map up to dimension $j$. 

\item The support of each simplex in the subdivision $\hat T\tau$ of $\tau$ is contained in some $U_k$.

\item Suppose that $\tau:\Delta^j\to X$ and that $\Delta'$ is the simplicial subdivision of $\tau$ such that $\hat T\tau$ is a singular subdivision of $\tau$ based on $\Delta'$. Suppose further that $\delta$ is a $j-1$ simplex of $\Delta'$ such that for some $l$, $0\leq l\leq j-1$,  an $l$-face $\beta$ of $\delta$ is contained in the $l$ skeleton of $\Delta^j$. Then, identifying $\beta$ with its singular subdivision inherited from that of $\tau$, $\psi(\tau\circ \beta)=\psi(\eta)$ for each simplex $\eta$ contained in the closed star of $\tau\circ\beta$ in $\hat T(\tau)$. In other words, every singular simplex $\eta$ in $\hat T\tau$ having $\tau\circ \beta$ as an $l$-face has the same value under $\psi$ as does $\tau\circ \beta$ itself. 

\end{enumerate}

We must now define $\hat T$ on each singular $i$-simplex $\sigma$. $\hat T$ has already been defined on $\bd \sigma$ on which it satisfies the induction hypotheses. We wish to show that we can define $\hat T$ on $\sigma$ so that $\hat T\sigma$ will also satisfy the stated properties. To do this, we need only take a sufficiently fine barycentric subdivision of $\sigma$ \emph{holding $\bd \sigma$ fixed}. This construction is discussed in Munkres \cite[\S 16]{MK} for simplicial complexes. Here we can apply the process to the singular case by singular subdivision of our model simplex $\Delta^i$: given the simplicial complex $K$ determined by the subdivision of $\bd \Delta^i$ induced by $\hat T(\bd \sigma)$, we first subdivide $\Delta^i$  compatibly with $K$ by taking $\bar c K$ to obtain $\Delta'$. Then we take a sufficiently iterated barycentric subdivision of $\Delta'$ mod $|K|=|\bd \Delta^i|$ as in \cite[\S 16]{MK}. This determines a singular subdivision of $\Delta^i$ by the partial ordering that preserves the existing partial orderings on the boundaries and then orders the new barycenters by the dimensions of the faces of which they are barycenters and by stage of construction (just as for ordinary iterated barycentric subdivision). We note that clearly condition 1 will hold, and  we will also have $\bd \hat T\sigma=\hat T\bd \sigma$, which provides condition 2. Conditions 3 and  4 can be achieved since these collectively impose a finite number of conditions on the degree of the subdivision that must be taken. In particular, it is not hard to see that  condition 3  can be satisfied by a direct application of \cite[Lemma 16.3]{MK}. For condition 4, we observe that, by induction, all singular simplices in the star of $\tau\circ\beta$ in $\hat T(\bd \tau)$ satisfy the given condition that they should evaluate to $\psi(\tau\circ\beta)$ under $\psi$. In particular, $\tau^{-1}(U_{\psi(\tau\circ\beta)})$ contains the simplicial star of $\beta$ in $K$, and no simplex in $K$ having  $\beta$ as a face is contained in any $U_k$ with $k<\psi(\tau(\beta))$ in the chosen ordering. The proof of \cite[Lemma 16.3]{MK} then demonstrates that under a sufficiently fine barycentric subdivision $L$ of $\Delta'$ mod $K$, the star of $\beta$ in $L$ will also be contained in $\tau^{-1}(U_{\psi(\tau\circ\beta)})$. This suffices to satisfy condition 4.

The preceding paragraph shows that we may obtain by induction a chain map $\hat T: C^c_*(X)\to {}_{\mc U}C^c_*(X)$. Furthermore, we can turn $\hat T$ into a map $T: I^{\bar p}C_*^{c}(X;\mc G_0 )\to I^{\bar p}C_*^{c}(X;\mc G_0 )$ as follows: Suppose that $\xi=\sum n_j \sigma_j$ and that $\hat T(\sigma_j)=\sum_l (-1)^{\text{sgn}(\i_{j,l})}\sigma_j \i_{j,l}$, where $\i_{j,l}$ are  singular simplices in the singular subdivision of $\Delta_j$. Then we set $T(\xi)=T(\sum n_j \sigma_j)=\sum_j \sum_l (-1)^{\text{sgn}(\i_{j,l})}(n_j\i_{j,l}) (\sigma_j\i_{j,l})$. Lemma \ref{L: sub a} shows that this map is well-defined on intersection chains since $T(\xi)$ is always a subdivision of $\xi$. 

We next need to show that each chain in the image of $T$ can be written as a sum of allowable chains each supported in some $U_j$. So let $\xi$ be an intersection $i$-chain. We will write $T\xi=\sum \xi_j$ with $|\xi_j|\subset U_j$. Since $\xi$ is a finite chain and $T$ takes only finite subdivisions of each simplex, $T\xi$ will also be a finite chain. We let $\xi_j$ be the subchain of $T\xi$ consisting of simplices (with coefficients) on which $\psi$ evaluates to $U_j$. We must show that each such $\xi_j$ is allowable.

If $\eta$ is an $i$-simplex in some $\xi_j$, then $\eta$ is a subdivision simplex of some $i$-simplex $\sigma$ of the allowable chain $\xi$. Hence $\eta$ is allowable by the arguments in Lemma \ref{L: sub a}. So it remains to consider the allowability of $i-1$ chains in $\bd \xi_j$.

Let $\mu$ be a singular $i-1$ simplex in $\bd \xi_j$. Then $\mu$ is an $i-1$ face of a singular $i$ simplex $\eta$, which is a singular $i$ simplex of $\hat T\sigma$ for some $\sigma:\Delta^i\to X$ in $\xi$. Furthermore $\psi(\eta)=j$. Let $\delta$ be the simplicial $i-1$ simplex corresponding to $\mu$ in the subdivision $\Delta'$ of $\Delta^i$ upon which $\hat T\sigma$ is based. First suppose that it is not true for any $l$, $0\leq l\leq i-1$, that $\delta$ has an $l$ face in the $l$ skeleton of $\Delta^i$. This implies that  $\delta$ intersects the $l$ skeleton of $\Delta^i$ only in its own $l-1$ skeleton. Thus  $\mu^{-1}(X_{n-k}-X_{n-k-1})=\delta\cap \sigma^{-1}(X_{n-k}-X_{n-k-1})\subset \delta \cap \{i-k+p(k) \text{ skeleton of }\Delta^i\} \subset \{ i-1+k+p(k) \text{ skeleton of }\delta\}$, and $\mu$ is allowable. 

Suppose next that there is some $l$, $0\leq l\leq i-1$, such that $\delta$ has an $l$ face $\beta$ in the $l$ skeleton of $\Delta^i$. Then by construction, every simplex in  the star of $\beta$ in $\Delta'$ gets taken under $\sigma$ into $U_{\psi(\sigma\circ\beta)}$ but not into any $U_m$ for $m<\psi(\sigma\circ\beta)$. In particular, since $\psi(\eta)=j$, $\psi$ evaluates to $j$ for all singular simplices built on simplices in the star of $\beta$ in $\Delta'$. This includes all simplices in the star of $\delta$, whence  
 the closed star of $\mu$  in $T\xi$ is also in $\xi_j$. So the coefficient of $\mu$ in $\bd \xi_j$ must be the same as that of $\mu$ in $\bd T\xi$.  But since $T\xi$ is allowable either $\mu$ is allowable or the coefficient of $\mu$ in $\bd T\xi$ is $0$, in which case $\mu$  must not be in $\bd \xi_j$. Either way, we see that $\bd \xi_j$ is allowable.

Thus we conclude that the image of $T$ is indeed in $I^{\bar p}_{\mc{U}}C^c_*(X;\mc {G})$. The desired chain homotopy $D$ from $\iota T$ to the identity can  be constructed as in the proof Proposition \ref{P: subdivision}. Since this time we have constructed our subdivision operator as a chain map (as opposed to our previous study of subdivisions simply on given chains), the inductive construction of $D$ in Proposition \ref{P: subdivision} provides a chain homotopy: we need only note that the  terms denoted $E$ in that proof can here be realized as $D(\bd \xi)$.

\end{proof}

\begin{corollary}\label{C: small equiv}
$\iota: I^{\bar p}C^c_*(X;\mc {G}_0) \to I^{\bar p}_{\mc{U}}C^c_*(X;\mc {G}_0)$ is a  chain homotopy equivalence, hence $I^{\bar p}H^c_*(X;\mc G_0)\cong H_*(I^{\bar p}_{\mc{U}}C^c_*(X;\mc {G}_0))$.
\end{corollary}
\begin{proof}
By Proposition \ref{P: U small}, there is a  chain homotopy $D$ from $\iota T$ to the identity. Consider then $T\iota$. $\iota$ is injective, being induced by inclusion, and $T\iota$  takes a chain $\xi\in I^{\bar p}_{\mc{U}}C^c_*(X;\mc {G}_0) $ and returns a subdivision. We also observe that the chain homotopy $D$ is well defined on the subcomplex $I^{\bar p}_{\mc{U}}C^c_*(X;\mc {G}_0) $ since for any allowable chain $\zeta$, $|D\zeta|\subset |\zeta|$. Thus we can define a chain homotopy $\bar  D$ on $I^{\bar p}_{\mc{U}}C^c_*(X;\mc {G}_0)$ by $\bar D=\iota^{-1}D\iota$ ($\iota^{-1}$ being well-defined on the subcomplex $I^{\bar p}_{\mc{U}}C^c_*(X;\mc {G}_0) $ by the injectivity of $\iota$). Then
$1-T\iota=\iota^{-1}\iota(1-T\iota)=\iota^{-1}(1-\iota T)\iota=\iota^{-1}(\bd D+D\bd)\iota=\bd  \bar D+\bar D\bd$. Thus $T$ and $\iota$ are chain homotopy inverses. 

\end{proof}

\begin{lemma}\label{L: excision2}
Let $X$ be a Hausdorff filtered space, let $U \subset X$  be any open subspace, and let $V$ be a closed subspace of $U$. Then $I^{\bar p}H^c_*(X,U;\mc G_0)\cong I^{\bar p}H^c_*(X-V, U-V;\mc G_0)$.
\end{lemma}
\begin{proof}
Using the preceding corollary, the proof now follows exactly as in the standard singular chain case; see, e.g., the proof of \cite[Theorem 31.7]{MK}.
\end{proof}

\subsection{Computations}\label{S: compute}

In this section we indulge in the computations that  make intersection homology theory go and which will enable us to perform the required verification that our theory satisfies the sheaf axioms if $X$ is a pseudomanifold. In particular, we here compute the intersection homology of products with $\R^n$, cones, distinguished  neighborhoods ($\cong cL\times \R^k$), and deleted distinguished  neighborhoods ($\cong (cL-x)\times \R^k$). Not surprisingly, the results presented here bear a marked similarity to those pre-existing in the literature (e.g. \cite{GM2,Bo, Ki}), however it is necessary that we proceed from scratch as these sources rely either on PL chains or compactly supported chains and, of course, they  assume traditional coefficient systems. We must proceed from first principles to derive these formulae for locally-finite singular chains with coefficients in $\mc G_0$. 

In most cases, our strategy will be to reduce our computation to that of intersection homology with compact supports and then proceed from there using the availability of stratum-preserving homotopy invariance in that setting.
The following lemma shows that the finite and locally-finite theories agree for intersection homology relative to a cocompact space. Of course the standard proof for ordinary homology would just involve breaking chains into a compact piece and a non-compact piece which can be thrown away. For intersection homology, however, we don't have such liberty to break chains (newly introduced boundaries may be in-allowable), but it turns out that we can break them in certain ways after performing a sufficient subdivision.

\begin{lemma}\label{L: cocompactness}
Let $X$ be a  filtered space with  coefficients $\mc G_0$, and let $U$ be an open subset such that $X-U$ is compact. Then  $I^{\bar p}H^{\infty}_*(X,U;\mc G_0)\cong I^{\bar p}H^c_*(X,U;\mc G_0)$.
\end{lemma}   
\begin{proof}
We begin with the obvious map $ I^{\bar p}C^c_*(X,U;\mc G_0)\to I^{\bar p}C^{\infty}_*(X,U;\mc G_0)$ induced by inclusion at the chain level and show that it induces a homology isomorphism.

We first show surjectivity: Let $\xi$ be an $i$-chain representing an element of $ I^{\bar p}H^{\infty}_i(X,U;\mc G_0)$. The ``obvious'' thing to do would be to cut out all of the simplices of $\xi$ with support in  $X-U$. However, this cannot be done directly, as the resulting boundaries may not be allowable (e.g., $\xi$ could be composed of an infinite number of  simplices with unallowable boundaries that just happen to cancel when taking the chain boundary). So we must refine the argument. 

 Consider the barycentric subdivision $\xi'$ of $\xi$. Since $\xi'$ and $\xi$ are relatively homologous by Proposition \ref{P: subdivision}, it suffices to find a finite chain relatively homologous to $\xi'$.  Let $\Xi$ denote the subchain consisting of  the singular simplices in $\xi$ (with their coefficients) whose supports intersect $X-U$. Note that $\Xi$ is comprised of a finite number of simplices since $X-U$ is compact. Let $\eta$ be the finite chain comprised of all singular $i$-simplices in $\xi'$ (with their coefficients) that share a vertex with a singular simplex of $\Xi$. Note that this includes all singular simplices in the barycentric subdivision $\Xi'$ of $\Xi$. We claim that $\eta$ and $\gamma=\xi'-\eta$ are each allowable chains. If so, then since $\eta$ contains all simplices from the subdivision of $\Xi$, $\gamma$ must have support in $U$, and so $\eta=\xi'$ in $ I^{\bar p}H^{\infty}_i(X,U;\mc G_0)$. 

To prove the claim, we first note that all simplices in $\xi'$ are allowable, as shown in the proof of Lemma \ref{L: sub a}. It remains to show that $\bd \eta$ is composed of allowable $i-1$ simplices, from which it will also follow that $\bd \gamma=\bd \xi'-\bd \eta$ is allowable, since $\bd \xi'$ is. The simplices in $\bd \eta$ all will be $i-1$ faces of $i$-simplices of $\xi'$. Up to orientation and vertex ordering convention, each singular $i$-simplex in $\xi'$ has the form of $\sigma\i$ with $\i: \Delta^i\to \Delta^i$ taking $\Delta^i$  linearly and injectively to  the polyhedral $i$-simplex $[\hat \Delta_i,\ldots, \hat \Delta_0]$, where $\hat \Delta_k$ is a the barycenter of a $k$-face of $\Delta$, the model simplex for $\sigma$ (we refer the reader to Munkres \cite{MK} for an exposition on barycentric subdivisions). Thus the $i-1$-simplices of $\bd \eta$ will similarly be compositions of $\sigma$ with singular $i-1$ simplices $\j:\Delta^{i-1}\to \Delta^i$ taking $\Delta^{i-1}$ to a polyhedral $i-1$ simplex $\delta\subset \Delta^i$ having  $i-1$ distinct vertices chosen from the set $\{\hat \Delta_i,\ldots, \hat \Delta_0\}$. We identify $\delta$ with the model simplex $\Delta^{i-1}$ via $\j$. 

Let us choose some such singular $i-1$ simplex $\tau=\sigma|_{\delta}$ (up to orientation and vertex ordering) in the boundary of a simplex of $\eta$. We consider separately the case of whether or not the simplex $\delta$ contains a vertex $\hat \Delta_0$. First, suppose not. In this case, $\delta$ has the form   $[\hat \Delta_i,\ldots, \hat \Delta_1]$ and the intersection of $\delta$ with the $m$ skeleton of $\Delta$ lies in the $m-1$ skeleton of $\delta$. This is because any simplex in the intersection of $\delta$ with the $m$ skeleton of $\Delta$ must be a face of a simplex  of the form $[\hat \Delta_m,\ldots, \hat \Delta_1]$, which is an $m-1$ simplex. Thus we see that $\tau^{-1}(X_{n-k}-X_{n-k-1})=\delta\cap \sigma^{-1}(X_{n-k}-X_{n-k-1})\subset \delta \cap \{i-k+p(k) \text{ skeleton of }\Delta\} \subset \{
 i-1+k+p(k) \text{ skeleton of }\delta\}$. Hence $\tau$ is allowable. 

Suppose on the other hand that $\delta$ includes a vertex $\hat \Delta_0$, which is a vertex of $\Delta$. Note that $\hat \Delta_0$ must be a vertex of a simplex in $\Xi$ or else $\tau$ would not be a simplex in $\eta$. Now since $\tau$ is an $i-1$ face of a simplex of $\xi'$ and since $\xi'$ is an allowable  relative cycle, we know that $\tau$ is either allowable or there are other singular $i$ simplices in $\xi'$ that also include $\tau$ as a boundary simplex and such that all the coefficients of $\tau$ in $\bd \xi'$ cancel (or else allowability of the boundary of $\xi'$ would be violated). In the first case (allowability) we are done. In the second case, we note that all of the other $i$ simplices that provide canceling boundary pieces are also in $\eta$, by our choice of $\eta$, since they will also have $\hat \Delta_0$ as a vertex, i.e. the full star of $\hat \Delta^0$ in $\xi'$ is in $\eta$. Thus the cancellation of $\tau$ also occurs in $\bd \eta$.  

Hence we have shown that all simplices in $\bd \eta$ are allowable, so $\eta$ and $\gamma$ are allowable and $\eta$ is a finite chain representing $\xi$ in $ I^{\bar p}H^{\infty}_i(X,U;\mc G_0)$. This proves surjectivity of the map on intersection homology induced by the inclusion $ I^{\bar p}C^c_*(X,U;\mc G_0)\to I^{\bar p}C^{\infty}_*(X,U;\mc G_0)$ .

For injectivity, suppose  two cycles $\xi_1, \xi_2\in I^{\bar p}C^c_*(X,U;\mc G_0)$ are relatively homologous via a chain $\zeta$ in $ I^{\bar p}C^{\infty}_*(X,U;\mc G_0)$. Then we apply the above procedure to first replace everything with subdivisions and then cut $\zeta'$ into two pieces $\eta$ and $\gamma$, the first finite and the second with support in $U$. Then $\bd \eta+\bd \gamma=\bd \zeta'=\xi_1'-\xi_2'+\omega$, where $\omega$ has support in $U$. So $\bd \eta-\xi_1'+\xi_2' =\omega -\bd \gamma $. Since $\bd \eta$, $\xi_1'$, and $\xi_2'$ are finite, $\omega -\bd \gamma$ must also be finite and its support is in $U$. Thus $\eta$ provides a relative homology from $\xi_1'$ to $\xi_2'$ in $I^{\bar p}C^c_*(X,U;\mc G_0)$. So the map $I^{\bar p}H^c_*(X,U;\mc G_0)\to I^{\bar p}H^{\infty}_*(X,U;\mc G_0)$ induced by inclusion is also injective.
\end{proof}

\begin{corollary}\label{C: compact rel}
Let $X$ be a compact filtered space with stratified coefficients $\mc G_0$.  Let $\R^*$ denote $(-\infty,0)\cup (0,\infty)\subset \R$.
Then  $I^{\bar p}H^{\infty}_*(X\times \R,X\times \R^*;\mc G_0\times \R)\cong I^{\bar p}H^c_*(X\times \R,X\times \R^*;\mc G_0\times \R)$.
\end{corollary}

In order to be able to apply the preceding lemma to compute locally-finite intersection homology of distinguished neighborhoods, we need to find a way to turn computations of absolute intersection homology groups into computations of intersection homology groups relative to cocompact subspaces. The following lemma is a first step towards making this possible by showing that certain cocompact subsets must have trivial intersection homology. By the long exact sequence of the pair, this will show that absolute and relative intersection homology agree for the cases of interest.

Recall that if $U$ is a subspace of $X$, then we define $I^{\bar p}C_*(U_X;\mc G_0)$ to be the chain subcomplex of $I^{\bar p}C_*(X;\mc G_0)$ consisting of allowable chains in $X$ with support in $U$.

\begin{lemma}\label{L: zero ends}
Let $L$ be a  filtered space with coefficients $\mc G_0$. Then $I^{\bar p}H^{\infty}_*((L\times \R^*)_{L\times \R};\mc G_0\times \R)=0$.
\end{lemma}
\begin{proof}
The lemma states that every allowable locally-finite cycle $\xi$ in $L\times \R$ with support in $L\times \R^*$ bounds an allowable locally-finite chain $\Xi$ in $L\times \R^*$. 

So let $\xi$ be such an $i$ cycle. We must construct $\Xi$. Suppose $\xi=\sum n_j\sigma_j$ with $\sigma_j:\Delta_j\to L\times \R^*$.  To construct $\Xi$, we begin with $\amalg \Delta_j$ and consider a locally finite singular triangulation of $\amalg \Delta_j\times [0,\infty)$  that gives the standard (non-subdivided)  singular triangulation of $\Delta_j\times 0$ with its orientation in $\bd(\Delta\times [0,\infty))$ and such that if $\sigma_l$ and $\sigma_k$  agree on an $i-1$ face  then corresponding triangulations on the products of those faces with $[0,\infty)$ will be compatible (we can build such a triangulation inductively in a standard way over the polyhedral skeleta of $\amalg \Delta_j$). Then if $\pi_L: L\times \R\to L$ and $\pi_{\R}:L\times \R\to \R$ are the projections, we can consider the map $f:\Delta_j\times \R \to L\times \R$ given by $(x,t)\to (\pi_L(\sigma_j(x)), \pi_{\R}(\sigma_j(x))\pm t)$, where the sign is $+$ if $\pi_{\R}(\sigma_j(x))>0$ and $-$ if $\pi_{\R}(\sigma_j(x))<0$. The map $f$ then determines the singular simplices of $\Xi$ by composition with the linear inclusions that give the singular $i+1$ simplices in the triangulation of $\amalg \Delta_j\times [0,\infty)$, and their coefficient lifts are determined from those of the original chain by the unique lifting property on $\mc G$. 

It should be clear that $\Xi$ provides the desired nullhomology of $\xi$ provided that $\Xi$ is allowable and locally-finite. The boundary of $\Xi$ is $\xi$, which we already know is allowable. Let $L\times \R=X$. The $i+1$ simplices of $\Xi$ are allowable since, if $\tau$ is such a simplex based upon the polyhedral  $i+1$ simplex $\delta\subset \Delta_j\times [0,\infty)$ (which we identify with the standard $i+1$ model simplex via its embedding in the singular subdivision), then  $\tau^{-1}(X_{n-k}-X_{n-k-1})=\delta\cap( \sigma^{-1}(X_{n-k}-X_{n-k-1})\times [0,\infty)) \subset \delta\cap [\{i-k+\bar p(k) \text{ skeleton of } \Delta_j\} \times [0,\infty) ]\subset \delta \cap \{i+1-k+\bar p(k) \text{ skeleton of } (\Delta_j\times [0,\infty))\}\subset  \{i+1-k+\bar p(k) \text{ skeleton of } \delta\}$. 

For the local-finiteness, suppose that $z\in L\times \R^*$ and that no neighborhood of $z$ intersects the supports of only a finite number of simplices of $\Xi$. Clearly $z\in L\times (-N,N)$ for some $N$; consider $Z=z\times [-N,N]$. The subspace $Z$ is compact and so it can be covered by a finite number of neighborhoods that intersect the supports of only a finite number of simplices of $\xi$. But this implies by basic topology that there is a tube of the form $W\times [-N,N]$, $W$ an open subset of $L$, such that $W\times [-N,N]$ intersects the supports of only a finite number of simplices in $\xi$. Now, we note that  if $\amalg_{k\in \mc K} \sigma_k$ is the finite subset of singular simplices in $\xi$ whose supports intersect $W\times [-N,N]$ then the only singular $i+1$ simplices of $\Xi$ whose supports can intersect $W\times [-N,N]$ are at most those defined via singular simplices of $\amalg \Delta_j\times [0,\infty)$ that intersect  $\amalg_{k\in \mc K}\Delta_k\times [0,N]$. This  is a finite collection, so $W\times (-N,N)$ is a neighborhood of $z$ that intersects the supports of only a finite number of simplices of $\Xi$, a contradiction. So the chain $\Xi$ must be locally-finite. 
\end{proof}

The next proposition computes the intersection homology of a product of a compact filtered space with $\R$. Without the more straightforward transversality results of PL theory, the prototype of this computation for PL spaces given in \cite[\S II]{Bo} does not readily carry over to the singular chain case. However, see the remark following the proof of the proposition.

\begin{proposition}\label{P: times R}
Let $L$ be a compact filtered space. Then $I^{\bar p}H^{\infty}_*(L\times \R;\mc G_0\times \R)\cong I^{\bar p}H_{*-1}(L;\mc G_0)$. 
\end{proposition}
\begin{proof}
The short exact sequence 
\begin{diagram}
0&\rTo& I^{\bar p}C^{\infty}_*((L\times \R^*)_{L\times \R};\mc G_0\times \R^*) &\rTo& I^{\bar p}C^{\infty}_*(L\times \R;\mc G_0\times \R) \\
&&{}&\rTo & I^{\bar p}C^{\infty}_*(L\times \R, L\times \R^*;\mc G_0\times \R) &\rTo &0
\end{diagram}
gives rise to a long exact sequence in intersection homology. By Lemma \ref{L: zero ends}, 
$I^{\bar p}H^{\infty}_*(L\times \R^*_{(L\times \R)};\mc G_0\times \R^*)=0$, so $I^{\bar p}H^{\infty}_*(L\times \R;\mc G_0\times \R)\cong I^{\bar p}H^{\infty}_*(L\times \R, L\times \R^*;\mc G_0\times \R)$. By  Lemma \ref{C: compact rel}, $I^{\bar p}H^{\infty}_*(L\times \R, L\times \R^*;\mc G_0\times \R)\cong I^{\bar p}H^{c}_*(L\times \R, L\times \R^*;\mc G_0\times \R)$. This allows us to finish the calculation using compact chains. 

Using the stratum-preserving homotopy invariance of compactly supported intersection homology, we see that $I^{\bar p}H^{c}_*(L\times \R;\mc G_0\times \R)\cong I^{\bar p}H^{c}_*(L;\mc G_0)$, and $I^{\bar p}H^{c}_*(L\times \R^*;\mc G_0\times \R^*)\cong \oplus_{i=1,2}I^{\bar p}H^{c}_*(L;\mc G_0)$. It is also clear from stratum-preserving homotopy equivalences  that the map induced by inclusion $I^{\bar p}H^{c}_*(L\times \R^*;\mc G_0\times \R^*) \to I^{\bar p}H^{c}_*(L\times \R;\mc G_0\times \R)$ is surjective and that the maps obtained by restricting to the summands are identical. Thus from the long exact sequence of compactly supported intersection homology, $I^{\bar p}H^c_{*}(L\times \R, L\times \R^*;\mc G_0\times \R)\cong I^{\bar p}H^c_{*-1}(L;\mc G_0)$, which suffices since $L$ is compact.
\end{proof}

\begin{remark}\label{R: chain times R}
It will be useful here, and after each of the following propositions, to keep track of how the intersection homology isomorphisms can be represented by chain maps. In this proposition, for example, suppose we have a cycle $\xi$ representing an element of $I^{\bar p}H_{*-1}(L;\mc G_0)$. A  chain in $I^{\bar p}H^{\infty}_*(L\times \R;\mc G_0\times \R)$ representing the image of the class $[\xi]$ under the isomorphisms of the proof is given by $\xi\times \R$, by which we mean the following: for each singular simplex $\sigma_j: \Delta_j\to L$, consider a triangulation of $\Delta_j\times \R$. We then replace this triangulation with a singular triangulation (see Section \ref{S: subdivision}). Define the singular chain $\sigma\times \R$ by $\sigma\times \R=\sum_k (-1)^{\text{sgn}(\i_k)}(\sigma_j\times \text{id}_{\R})\circ \i_k$, where $\i_k$ are the simplices of the singular triangulation of $\Delta_j\times \R$ and $\sigma_j\times \text{id}_{\R}: \Delta_j\times \R\to L\times \R$ is  the product mapping. 
Choosing compatible triangulations for all $\Delta_j$ in $\xi$ (this can be done inductively by a uniform procedure at each dimension) and making the obvious corresponding modifications on coefficients gives a chain $\xi\times \R=\sum (n_j\times \R)(\sigma_j\times \R)$, which is locally-finite. It is also easily verified that $\xi\times \R $ is allowable. 

To see that $\xi\times \R$  indeed represents the image of $[\xi]$ in $I^{\bar p}H^{\infty}_*(L\times \R;\mc G_0\times \R)$, the important point to note is that for almost every $a\in (0,\infty)$, we can cut our  triangulation of $\Delta_j\times \R$ transversely at  $(\sigma_j\times \text{id}_{\R})^{-1}(\{-a,a\})$. So  further compatible subdivisions (simplicial then singular) of $\amalg \Delta_j\times [-a,a]$ followed by restriction of $\sigma_j\times \text{id}_{\R}$ yields a chain $\xi\times [-a,a]$ whose boundaries are subdivisions of $\xi \times a$ and $-\xi \times -a$.   Allowability of this new chain follows from that of $\xi$ and from the construction. That this is the correct chain follows from tracing through the isomorphisms of the proof and recalling that subdivision does not change intersection homology class (by Proposition \ref{P: subdivision}).

An alternative proof of  Proposition \ref{P: times R} would begin with such a map $\xi \to \xi\times \R$ and show directly that it induces an isomorphism. In fact, this is usually  what is done in the PL case (e.g. \cite[Ch. II]{Bo}), but attempts to mimic such proofs for singular chains encounter some difficulty. However, given a posteriori the isomorphisms of our proof, it is not difficult to perform the reverse engineering that gives us such a chain correspondence.\qedsymbol

\end{remark}

We must also compute the intersection homology of cones.
We think of $cL$ as $L\times [0,1)/ (x,0)\sim (y,0)$. The cone point is taken as the $0$ skeleton $(cL)_0$, and for $k>0$, $(cL)_k=L_{k-1}\times (0,1)$. We denote the induced stratified local coefficient system on $cL$ by $c\mc G_0$. This coefficient system is $\mc G\times (0,1)$ on $(L-\Sigma)\times (0,1)$ and $0$ elsewhere. Also, we will use $\bar cZ$ to denote the closed cone on $Z$: $\bar cZ=(Z\times [0,1])/(x,0)\sim (y,0)$.

\begin{remark}
We should note that  the following, seemingly innocuous, computation includes the crucial use of our two-tiered coefficient system, marking a deviation from traditional intersection homology computations. The main point is that $0$-cycles  behave very differently under coning than do higher dimensional cycles, for coning a point creates a $1$-chain with a new boundary component at the cone point. For traditional perversities, a cone on a $0$-cycle will never be allowable, which works out compatibly with the axioms for intersection homology. For superperversities, however, satisfaction of the intersection homology axioms requires that cones on $0$-cycles must be allowable, and in order for that to happen, the cone point boundary must vanish.  
\end{remark}

\begin{proposition}\label{P: cone}
Let $L$ be an $n-1$ dimensional filtered space with coefficient system $\mc G_0$. Then \begin{equation*}
I^{\bar p}H^{c}_i(cL;c\mc G_0) \cong
\begin{cases}
I^{\bar p}H^{c}_{i}(L;\mc G_0), & i<n-1-\bar p(n)\\
0, & i\geq n-1-\bar p(n).
\end{cases}
\end{equation*}
If $L$ is compact, then 
\begin{equation*}
I^{\bar p}H^{\infty}_i(cL;c\mc G_0) \cong
\begin{cases}
I^{\bar p}H_{i-1}(L;\mc G_0), & i\geq n-\bar p(n)\\
0, & i<n-\bar p(n).
\end{cases}
\end{equation*}

\end{proposition}
\begin{proof}

To compute $I^{\bar p}H^c_*(cL;c\mc G_0)$, we argue as in \cite{Ki} and begin by determining which chains can intersect the $0$ stratum $(cL)_0$, which is the cone point, $x$. An allowable simplex $\sigma:\Delta^i\to cL$ must satisfy the condition that $\sigma^{-1}(x)$ be contained in the $i-n+\bar p(n)$ skeleton of $\Delta^i$. 

We first show that $IH_{i-1}^c(cL;c\mc G_0)=0$ if $i-n+\bar p(n)\geq 0$. If $\sigma:\Delta^{i-1}\to cL$ is a singular $i-1$ simplex, we define $\bar c \sigma:\bar c \Delta^{i-1}\cong \Delta^i\to cL$ as the singular simplex that takes the ray in $\bar c \Delta^{i-1}$ running from the cone point to $z\in \Delta^{i-1}$ linearly onto the ray in $cL$ running from the cone point to $\sigma(z)$. If $n$ is a coefficient of $\sigma$, the coefficient $\bar cn$ of $\bar c \sigma$ is determined by the homotopy lifting property on $\mc G\times (0,1)\subset c\mc G_0$. If $\xi=\sum n_j\sigma_j\in  IC_{i-1}^c(cL;c\mc G_0)$, we can then define $\bar c \xi\in IC_{i}^c(cL;c\mc G_0)$ by $\bar c\xi=\sum (\bar cn_j)(\bar c\sigma_j)$. 

Now, we have $i-n+\bar p(n)\geq 0$ if and only if $i\geq n-\bar p(n)$. So if $\xi$ is an allowable $i-1$ cycle with $i\geq n-\bar p(n)$, $i-1\geq 0$, then the $i$ chain $\bar c\xi$ will also be allowable: 
\begin{itemize}
\item its boundary is the allowable chain $\xi$ (even if $\xi$ is a $0$ cycle!)

\item if $k<n$,  each $i$ simplex $\bar c\sigma_j$ in $\bar c\xi$ satisfies 
\begin{align*}
(\bar c\sigma_j^{-1})((cL)_{n-k}-(cL)_{n-k-1})&= \sigma_j^{-1}((cL)_{n-k}-(cL)_{n-k-1})\times (0,1]\\
&\subset \{i-1-k+\bar p(k) \text{ skeleton of }\Delta^{i-1}_j\}\times (0,1]\\
&\subset  \{i-k+\bar p(k) \text{ skeleton of }\bar c\Delta^{i-1}_j\} 
\end{align*}

\item for $k=n$ (the stratum $(cL)_0=x$), we see that 
\begin{align*}
(\bar c\sigma_j)^{-1}(x)&=\bar c(\sigma_j^{-1}(x))\\
&\subset \bar c\{i-1-n+\bar p(n) \text{ skeleton of }\Delta^{i-1}_j\}\\
&\subset \{i-n+\bar p(n)  \text{ skeleton of }\bar c\Delta^{i-1}_j\cong \Delta^i \}.
\end{align*}

\end{itemize}

If $i-n+\bar p(n)<0$, then neither the $i-1$ cycle $\xi$ nor any potential $i$ chain whose boundary is $\xi$ can intersect the cone point $x$. So in this range $I^{\bar p}H_{i-1}^c(cL;c\mc G_0)\cong I^{\bar p}H^c_{i-1}(L\times (0,1);c\mc G_0)$ and stratum-preserving homotopy equivalence tells us that $I^{\bar p}H^c_{i-1}(cL;c\mc G_0 )\cong I^{\bar p}H^c_{i-1}(L;\mc G_0)$. 

This finishes the calculation if $I^{\bar p}H^c_*(cL;c\mc G_0)$.

To compute $I^{\bar p}H^{\infty}_*(cL;c\mc G_0)$, we first observe that  $I^{\bar p}H^{\infty}_*(cL;c\mc G_0)\cong I^{\bar p}H^{\infty}_*(cL, L\times (0,1);c\mc G_0)$. This will follow from the long exact sequence of the pair if  $I^{\bar p}H^{\infty}_*((L\times (0,1))_{cL};c\mc G_0)=0$, but  this can be seen  just as in the proof of Lemma \ref{L: zero ends} by ``pushing chains to infinity''. 

Next, by Lemma \ref{L: cocompactness}, we have an isomorphism $I^{\bar p}H^{\infty}_*(cL, L\times (0,1);c\mc G_0)\cong I^{\bar p}H^{c}_*(cL, L\times (0,1);c\mc G_0)$. We  compute $I^{\bar p}H^{c}_*(cL, L\times (0,1);c\mc G_0)$ via the long exact sequence of the pair for compact intersection chains. By stratum-preserving homotopy equivalence (Lemma \ref{L: sphe}), $I^{\bar p}H^{c}_*( L\times (0,1);c\mc G_0)\cong I^{\bar p}H^{c}_*(L;\mc G_0)$. Since $I^{\bar p}H^c_i(cL; c\mc G_0)=0$ if $i\geq n-1-\bar p(n)$, we see that $I^{\bar p}H^c_{i+1}(cL,L\times(0,1);c\mc G_0)\cong I^{\bar p}H^c_i(L\times(0,1);c\mc G_0)\cong  I^{\bar p}H^c_i(L;\mc G_0 )$ in this range. For  $i<n-1-\bar p(n)$, the map induced by inclusion $I^{\bar p}H^c_i(L;\mc G_0)\to I^{\bar p}H^c_i(cL;c\mc G_0)$ is an isomorphism, so $I^{\bar p}H^c_{i}(cL,L\times(0,1);c\mc G_0)=0$ for $i\leq n-1-\bar p(n)$. 

 The proposition now follows from these calculations. 
\end{proof}

\begin{remark}
Again, we would like an explicit chain construction of the isomorphism $I^{\bar p}H_i(cL;c\mc G_0)\cong I^{\bar p}H_{i-1}(L;\mc G_0)$ for $i\geq n-\bar p(n)$, $L$ compact.  So let $\xi=\sum_j n_j\sigma_j$ be a chain representing an element of $I^{\bar p}H_{i-1}(L;\mc G_0)$, $i\geq n-\bar p(n)$. Since $L$ is compact, $\xi$ will be a finite chain. This time for each $\sigma_j:\Delta_j\to L$, we consider $c\sigma_j: c\Delta_j \to cL$. Recall that $c$  denotes an open cone,  so $c\Delta_j\cong \Delta_j\times [0,1)/(x,0)\sim (y,0)$. Choosing a singular triangulation of $c\Delta_j$ and composing with the map $c\Delta_j \to cL$ that takes $(x,t)$ to $(\sigma_j(x), t)$ defines the chain $c\sigma_j$. Choosing compatible triangulations of each $c\Delta_j$ in $\xi$ and treating the coefficient lifts similarly allows us to define  a map $\xi\to c\xi$, which gives a chain representative of the image of $[\xi]$ in  $I^{\bar p }H_i^{\infty}(cL;\mc G_0)$ under the isomorphisms of the proof of the proposition. Once again (see Remark  \ref{R: chain times R}), this can be seen by making an appropriate transverse cut. 
\end{remark}

Putting together the previous calculations, we can compute the intersection homology of spaces of the form $cL\times \R^k$, the homeomorphism type of distinguished neighborhoods in pseudomanifolds.

\begin{proposition}\label{P: dist ngbd}
Let $L$ be a compact filtered space with coefficients $\mc G_0$. Then $I^{\bar p}H^{\infty}_*(cL\times \R^k;c\mc G_0\times \R^k)\cong I^{\bar p}H^{\infty}_{*-k}(cL;c\mc G_0)$. 
\end{proposition}
\begin{proof}
Let $x$ be the cone point of $cL$. For convenience, in this proof we treat $cL$ as $L\times [0,\infty)/(y,0)\sim(z,0)$. 

We begin with the claim that $I^{\bar p}H^{\infty}_*(cL\times \R^k;c\mc G_0\times \R^k)\cong I^{\bar p}H^{\infty}_*(cL\times \R^k, cL\times \R^k-(x,0);c\mc G_0\times \R^k)$. This will follow from the long exact sequence of the pair once we show that 
$I^{\bar p}H^{\infty}_*([cL\times \R^k-(x,0)]_{cL\times \R^k};c\mc G_0\times \R^k)=0$.  The proof of this fact follows from the same concepts as used in Lemma \ref{L: zero ends}. The principal tool in creating for each cycle $\xi=\sum n_j\sigma_j\in I^{\bar p}C^{\infty}_i([cL\times \R^k-(x,0)]_{cL\times \R^k};c\mc G_0\times \R^k)$  a chain $\Xi$ with $\bd \Xi=\xi$ is the use of a proper map $f:\amalg \Delta_j\times [0,\infty)\to cL\times \R^k-(x,0)$ built upon the original singular simplices $\sigma_j:\Delta_j \to cL\times \R^k-(x,0)$. Let us label points in $cL\times \R^k$ as $(z,s,v)$, where $z\in L$, $s\in [0,\infty)$, and $v\in \R^k$. Then for $y\in \Delta_j$, if $\sigma_j(y)=(z_0,s_0,v_0)$, we let $f(y, t)= (z_0, (1+t)s_0, (1+t)v_0)$. So $f(y,0)=\sigma_j(y)$, and as $t$ goes to $\infty$, $f(y,t)$ goes properly to the end of the space since $s_0$ and $v_0$ cannot both be $0$. Note  that $f$ is a stratum-preserving open-ended homotopy from $\amalg \sigma_j$. Thus arguments similar to those in Lemma \ref{L: zero ends} show that we can build an allowable $\Xi$ with $\bd \Xi=\xi$. Local-finiteness of $\Xi$ also follows  using the fact that $cL\times \R^k$ can be built as the increasing union of compact sets of the form $$\left\{\left(L\times [0,  \frac{2\tan^{-1}(N)}{\pi}]\right)/(y,0)\sim(z,0)\right\}\times [-N,N]^k.$$

So $I^{\bar p}H^{\infty}_*(cL\times \R^k;c\mc G_0\times \R^k)$ is isomorphic to $I^{\bar p}H^{\infty}_*(cL\times \R^k, cL\times \R^k-(x,0);c\mc G_0\times \R^k)$, which by Lemma \ref{L: cocompactness} is isomorphic to 
$I^{\bar p}H^{c}_*(cL\times \R^k, cL\times \R^k-(x,0);c\mc G_0\times \R^k)$. 

At this point, the standard way to proceed would utilize a K\"unneth Theorem or a Mayer-Vietoris sequence. These tools exist (as, in fact, will follow from our proof that we are indeed working with the  intersection homology modules as given  via sheaf theory), but rather than develop these tools, which would require some work, we instead proceed using induction on established results. 

As the base step, suppose $k=1$, and consider the short exact chain sequence (with $\R^*=\R-0$ and suppressing coefficients for readability)

\begin{diagram}
0 & \rTo & I^{\bar p}C^{c}_*(cL\times \R^*, (cL-x)\times \R^*)
& \rTo & I^{\bar p}C^{c}_*(cL\times \R, (cL-x) \times \R)\\
{}& \rTo & I^{\bar p}C^{c}_*(cL\times \R)/\{I^{\bar p}C^{c}_*((cL-x) \times \R)+ I^{\bar p}C^{c}_*(cL\times \R^*)\}
& \rTo & 0,&	
\end{diagram}
where we have replaced the kernel $$\{I^{\bar p}C^{c}_*((cL-x) \times \R)+ I^{\bar p}C^{c}_*(cL\times \R^*)\}/I^{\bar p}C^{c}_*((cL-x) \times \R)$$ with the isomorphic 
\begin{align*}
I^{\bar p}C^{c}_*(cL\times \R^*)/\{I^{\bar p}C^{c}_*(cL\times \R^*)\cap I^{\bar p}C^{c}_*((cL-x) \times \R)\}&\cong I^{\bar p}C^{c}_*(cL\times \R^*)/I^{\bar p}C^{c}_*( (cL-x)\times \R^*)\\
&=I^{\bar p}C^{c}_*(cL\times \R^*, (cL-x)\times \R^*).
\end{align*}

This sequence yields a long exact sequence in homology. It follows from Proposition \ref{P: U small} that the inclusion $\iota: I^{\bar p}C^{c}_*((cL-x) \times \R)+ I^{\bar p}C^{c}_*(cL\times \R^*)\into I^{\bar p}H^{c}_*(cL\times \R^k-(x,0))$ is a chain homotopy equivalence. Thus the homology of the quotient  term in the short exact sequence is isomorphic to $I^{\bar p}H^c_*(cL\times \R,(cL\times \R)-(x,0))$, and the associated long exact sequence in homology is  

\begin{diagram}
0 & \rTo & I^{\bar p}H^{c}_*(cL\times \R^*, (cL-x)\times \R^*)
& \rTo & I^{\bar p}H^{c}_*(cL\times \R, (cL-x) \times \R)\\
{}& \rTo & I^{\bar p}H^c_*(cL\times \R,(cL\times \R)-(x,0))
& \rTo & 0&,	
\end{diagram}

By stratum preserving homotopy invariance of compactly supported intersection homology (Lemma \ref{L: sphe}), 
$$ I^{\bar p}H^{c}_*(cL\times \R, (cL-x) \times \R;c\mc G_0\times \R)\cong I^{\bar p}H^{c}_*(cL, cL-x;c\mc G_0),$$ $$I^{\bar p}H^{c}_*(cL\times \R^*, (cL-x)\times \R^*;c\mc G_0\times \R)\cong \oplus_{i=1,2}I^{\bar p}H^{c}_*(cL, cL-x;c\mc G_0),$$ 
and the map from the latter to the former is an isomorphism  on restriction to each summand. So the long exact sequence splits into split short exact sequences, which shows that $$I^{\bar p}H^{c}_*(cL\times \R, (cL\times \R)-(x,0);c\mc G_0\times \R)\cong I^{\bar p}H^{c}_{*-1}(cL, cL-x;c\mc G_0).$$ The proposition now follows for $k=1$ using $$I^{\bar p}H^{c}_{*-1}(cL, cL-x;c\mc G_0)\cong I^{\bar p}H^{\infty}_{*-1}(cL, cL-x;c\mc G_0)\cong I^{\bar p}H^{\infty}_{*-1}(cL;c\mc G_0)$$ 
(see the proof of Proposition \ref{P: cone}).

Now suppose inductively that 
$$I^{\bar p}H^{c}_*(cL\times \R^{j-1}, cL\times \R^{j-1}-(x,0))\cong I^{\bar p}H^{c}_{*-(j-1)}(cL, cL-x)$$ for $j<k$. Consider now the long exact sequence associated to the short exact sequence
\begin{diagram}
0 & \rTo & I^{\bar p}C^{c}_*(cL\times \R^{k-1}\times  \R^*, (cL\times \R^{k-1}-(x,0_{k-1}) )\times \R^*)\\
{}& \rTo & I^{\bar p}C^{c}_*(cL\times \R^{k-1}\times \R, (cL\times \R^{k-1}-(x,0_{k-1}) ) \times \R)\\
{}& \rTo & I^{\bar p}C^{c}_*(cL\times \R^k)/\{I^{\bar p}C^{c}_*((cL\times \R^{k-1}-(x,0_{k-1}) ) \times \R) + I^{\bar p}C^{c}_*(cL\times \R^{k-1}\times  \R^*)\}
& \rTo & 0.	
\end{diagram}

Once again we have used the obvious isomorphisms to write the first term in a convenient form, and, using Proposition \ref{P: U small}, the homology of the third term is just $I^{\bar p}H^c_*(cL\times \R^k,cL\times \R^k-(x,0))$. 
Also as before, stratum preserving homotopy equivalence and an exact sequence argument tell us that 
$$I^{\bar p}H^{c}_*(cL\times \R^{k-1}\times \R, (cL\times \R^{k-1}-(x,0_{k-1}) ) \times \R;) \cong I^{\bar p}H^{c}_*(cL\times \R^{k-1}, (cL\times \R^{k-1}-(x,0_{k-1}) ) ), $$
$$I^{\bar p}H^{c}_*(cL\times \R^{k-1}\times  \R^*, (cL\times \R^{k-1}-(x,0_{k-1}) )\times \R^*)\cong \oplus_{i=1,2}I^{\bar p}H^{c}_*(cL\times \R^{k-1}, (cL\times \R^{k-1}-(x,0_{k-1}) ) ), $$
and then 
$$I^{\bar p}H^{c}_*(cL\times \R^k, cL\times \R-(x,0))
\cong I^{\bar p}H^{c}_{*-1}(cL\times \R^{k-1}, (cL\times \R^{k-1}-(x,0_{k-1}) ) ).$$ Applying the induction hypothesis, we see that $$I^{\bar p}H^{c}_*(cL\times \R^{k}, cL\times \R^{k}-(x,0))\cong I^{\bar p}H^{c}_{*-k}(cL, cL-x),$$ and the rest of the theorem follows as in the case $k=1$. 

\end{proof}

\begin{remark}
If $\xi$ is a chain representing $[\xi]\in I^{\bar p}H^{\infty}_{*-k}(cL;c\mc G_0)$ then the image of $[\xi]$ in  $I^{\bar p}H^{\infty}_*(cL\times \R^k;c\mc G_0\times \R^k)$ can be represented by $\xi\times \R^k$, where $\xi\times \R^k$ is constructed analogously to  $\xi\times \R$ in Remark \ref{R: chain times R} via triangulation of  each $\Delta_j\times \R^k$. \qedsymbol
\end{remark}

We also want to compute the intersection homology of deleted distinguished neighborhoods:

\begin{proposition}\label{P: hole}
Let $L$ be a compact filtered space with coefficients $\mc G_0$. Then $$I^{\bar p}H^{\infty}_*((cL-x)\times \R^k;c\mc G_0\times \R^k)\cong I^{\bar p}H^{\infty}_{*-k}(cL-x;c\mc G_0).$$
\end{proposition}
\begin{proof}
We note that $cL-x\cong L\times \R$ and $(cL-x)\times \R^k\cong L\times \R^{k+1}$. To see that the appropriate intersection homology groups are isomorphic, we proceed just as in the proceeding proposition. In particular, 
\begin{align*}I^{\bar p}H^{\infty}_*((cL-x)\times \R^k;c\mc G_0\times \R^k)&\cong 
I^{\bar p}H^{\infty}_*(L\times \R^{k+1};c\mc G_0\times \R^k)\\
&\cong I^{\bar p}H^{\infty}_*(L\times \R^{k+1}, L\times (\R^{k+1}-0);c\mc G_0\times \R^k)\\
&\cong I^{\bar p}H^{c}_*(L\times \R^{k+1}, L\times (\R^{k+1}-0);c\mc G_0\times \R^k).
\end{align*}
Then again we induct, this time using the short exact sequence
\begin{diagram}
0 & \rTo & I^{\bar p}C^{c}_*(L\times \R^{k}\times  \R^*, L\times (\R^{k}-0_{k})\times \R^*)\\
{}& \rTo& I^{\bar p}C^{c}_*(L\times \R^{k}\times \R, L\times (\R^{k}-0_{k}) \times \R)\\
{}& \rTo & I^{\bar p}C^{c}_*(L\times \R^{k+1})/\{ I^{\bar p}C^{c}_*( L\times (\R^{k}-0_{k})  \times \R)+ I^{\bar p}C^{c}_*(L\times \R^{k}\times  \R^*)\}
& \rTo & 0.&	
\end{diagram}
\end{proof}

\begin{remark}
In chains, we can again take $\xi\in I^{\bar p}C^{\infty}_{*-k}(cL-x)$ to $\xi\times \R^{k}\in I^{\bar p}C^{\infty}_{*}((cL-x)\times \R^k)$. 
\end{remark}

The final calculation of this section, contained in the following lemma and corollary, establishes that local intersection homology of a pseudomanifold can be computed through the use of a single distinguished neighborhood. In other words, we show that $\lim_{x\in U}I^{\bar p}H^{\infty}_*(X,X-\bar U;\mc G_0)$ is the direct limit of an essentially constant direct system with a cofinal set consisting of distinguished neighborhoods, the maps between which are intersection homology isomorphisms.

\begin{lemma}\label{L: rest dist}
Let $X$ be a pseudomanifold, $x\in X$, and let $N$ be a distinguished neighborhood of $x$, i.e. $N\cong cL\times \R^{n-k}$. Assume also that $\bar N$ is compact and homeomorphic to $\bar cL\times D^{n-k}$, where $\bar cL$ is the closed cone on $L$ and $D^{n-k}$ is the closed unit disk in $\R^{n-k}$. For $\alpha\in (0,1)$, let $N_{\alpha}\subset N$ be a distinguished neighborhood of $x$ in $N$ such that if $\phi: \bar cL\times D^{n-k}\to \bar N$ is the homeomorphism, then $N_{\alpha}=\phi(\alpha cL\times \alpha D^{n-k})$, where $\alpha D^{n-k}$ is the open subdisk of $D^{n-k}$ of radius $\alpha$ and $\alpha cL=L\times [0,\alpha)/(y,0)\sim(z,0)\subset cL=L\times [0,1]/(y,0)\sim(z,0)$. Then if $\alpha<\beta\in (0,1)$, the natural quotient $I^{\bar p}C^{\infty}_*(X,X-\bar N_{\beta};\mc G_0)\to I^{\bar p}C^{\infty}_*(X,X-\bar N_{\alpha};\mc G_0)$ induces an isomorphism on intersection homology.
\end{lemma}
\begin{proof}
We consider the exact sequence of the triple $(X,X-\bar N_{\alpha}, X-\bar N_{\beta})$. Then it suffices to show that $I^{\bar p}H^{\infty}_*((X-\bar N_{\alpha})_X, (X-\bar N_{\beta})_X;\mc G_0)=0$. Note that all chains must be locally-finite in $X$. We first claim that $I^{\bar p}H^{\infty}_*((X-\bar N_{\alpha})_X, (X-\bar N_{\beta})_X;\mc G_0)\cong I^{\bar p}H^{c}_*(X-\bar N_{\alpha}, X-\bar N_{\beta};\mc G_0)$. This follows  as in the proof of Lemma \ref{L: cocompactness} - since any chain in $I^{\bar p}C_*^{\infty}((X-\bar N_{\alpha})_X;\mc G_0)$ must be allowable in $X$, it must possess only a finite number of simplices with non-zero coefficients and supports intersecting $\bar N_{\beta}$. We can then proceed as in  Lemma \ref{L: cocompactness} to subdivide and truncate off a cofinite number of simplices supported in $X-\bar N_{\beta}$. But now $I^{\bar p}H^{c}_*(X-\bar N_{\alpha}, X-\bar N_{\beta};\mc G_0)=0$ as the two sets are stratum-preserving homotopy equivalent. 
\end{proof}

\begin{corollary}\label{C: essentially constant}
Suppose $X$ is a pseudomanifold, and $x\in X$. The direct system $I^{\bar p}H^{\infty}_*(X,X-\bar U;\mc G_0)\to I^{\bar p}H^{\infty}_*(X,X-\bar V;\mc G_0)$ determined by open sets $x\in V\subset U$ is essentially constant. In particular, for all $x\in X$, there exists a neighborhood $W\ni x$ such that $I^{\bar p}H^{\infty}_*(X,X-\bar W;\mc G_0)\cong \lim_{x\in U}I^{\bar p}H^{\infty}_*(X,X-\bar U;\mc G_0)$.
\end{corollary}
\begin{proof}
Let $N$ be a distinguished neighborhood of $x\in X$ as described in Lemma \ref{L: rest dist}; all points have such a neighborhood by shrinking $N$ if necessary. Then the neighborhoods $N_{\alpha}$ are cofinal, and the lemma states that the restriction map on such neighborhoods induces an isomorphism on homology. Now take $W=N_{\alpha}$ for any $\alpha\in(0,1)$. 
\end{proof}

\section{Sheaves} \label{S: sheaves}

In this section of the paper we construct a differential graded complex of sheaves based upon our singular intersection chain complex.  Although this complex will not satisfy the strongest properties we might want (softness, flabbiness, injectivy, etc.), it will be a homotopically fine sheaf complex, which will suffice to show that its hypercohomology  agrees (up to a reindexing) with the intersection homology modules we have already studied. Our main result is that on a paracompact stratified topological pseudomanifold, this sheaf complex is quasi-isomorphic to the Deligne sheaf complex, and hence its hypercohomology also agrees with Goresky-MacPherson sheaf intersection homology. In particular, if $\bar p$ is a superperversity, we obtain the intersection homology modules occuring in the superduality theorem of Cappell and Shaneson \cite{CS}.

\subsection{Definition and basic properties}\label{S: sheafify}

We fix a filtered \emph{Hausdorff} space $X=X^n\supset X^{n-1}\supset \cdots \supset X^0\supset X^{-1}=\emptyset$, a perversity or superperversity, and a  coefficient system $\mc G_0$, but we will omit these from the notation where there will be no confusion. 

We will consider two differential graded presheaves with cohomological indexing: $IS^*$ and $KS^*$. We define $IS^*$ by $U\to IC^{\infty}_{n-*}(X, X-\bar U)$ and $KS^*$ by $U\to IC^{c}_{n-*}(X, X-\bar U)$. The restriction maps are the obvious quotients in both cases, and there is a natural inclusion-induced presheaf morphism $i:KS^*\to IS^*$. These presheaves give rise to sheaves $\mc{KS^*}$ and $\mc{IS^*}$ with an induced morphism $i: \mc{KS^*}\to \mc{IS^*}$. In fact, this is a sheaf isomorphism:

\begin{lemma}\label{L: same sheaf}
The homomorphism $i: \mc{KS^*}\to \mc{IS^*}$ is an  isomorphism of sheaves. 
\end{lemma}
\begin{proof}
We must show that $i$ induces an isomorphism at each stalk. 

First, we show injectivity. Let $x\in X$ and $s\in \mc{KS}^{n-j}_x$, the stalk at $x$. Suppose that $U$ is a neighborhood of $x$ and  $\xi\in IC^c_j(X, X-\bar U;\mc G_0)$ is a finite  chain that represents $s$. If $i|_x(s)=0$, then $\xi=0$ in $IC^{\infty}_j(X,X-\bar V;\mc G_0)$ for some open $V$ such that $x\in V\subset U$. But this would imply that $|\xi|\in X-\bar V$, which implies that $\xi=0$ in $ IC^c_j(X, X-\bar V;\mc G_0)$. Hence $s=0$. 

For surjectivity, let $s\in \mc{IS}^{n-j}_x$, and suppose $U$ is a neighborhood of $x$ and $\xi\in IC^{\infty}_j(X,X-\bar U;\mc G_0)$ represents $s$. By taking a smaller $U$ if necessary, we may assume that $\bar U$ intersects the supports of only a finite number of the simplices of $\xi$. It will suffice to find a finite chain $\zeta\in IC^c_j(X,X-\bar U;\mc G_0)$ such that $i(\zeta)=\xi \in IC^{\infty}_j(X,X-\bar U;\mc G_0)$. If $\xi$ is already a finite chain then $\zeta=\xi$ suffices. Suppose then that $\xi$ contains an infinite number of singular simplices. Let $\Xi$ be the singular chain (not necessarily allowable) composed of singular simplices of $\xi$ (with their coefficients) whose supports intersects $\bar U$. Let $\xi'$ be the generalized barycentric subdivision of $\xi$ \emph{holding $\Xi$ fixed}. In other words, we perform  a barycentric subdivision of each simplex in $\xi$ except that we do not subdivide the simplices of $\Xi$ nor any common faces between simplices  in $\Xi$ and simplices not in $\Xi$ (see  \cite[\S 16]{MK}). Now take as $\zeta$ the ``regular neighborhood'' of $\Xi$ in $\xi'$. By this we mean take the chain consisting of the simplices in $\Xi$ (with their coefficents) and all other simplices in $\xi'$ that share a vertex  with a simplex in $\Xi$. This $\zeta$ must be finite since $\xi$ is locally finite and $\bar U$ is compact. Furthermore, $\zeta$ is allowable by exactly the same arguments as in the proof of Lemma \ref{L: sub a}. To see that $i(\zeta)=\xi$ in $IC^{\infty}_j(X,X-\bar U;\mc G_0)$, we simply note that $\xi-i(\zeta)$ has support in $|\xi-\Xi|$, which lies in $X-\bar U$. Hence $\xi-i(\zeta)=0$ in $IC^{\infty}_j(X,X-\bar U;\mc G_0)$.
\end{proof}

We would like to be able to say that the global sections of the presheaf $IS^*$ and of the sheaf $\mc{IS^*}$ agree, i.e. $IS^*(X)\cong \Gamma(X,\mc{IS}^*)$.  For such a statement to hold, it is only necessary that $IS^*$ have no non-trivial global sections with empty support and that it be conjunctive with respect to coverings (see \cite[I.6.2]{Br}). This is the content of the following lemmas.

Note that if $X$ is not compact, we do not expect  $KS^*(X)\cong \Gamma(X,\mc{IS}^*)$ since global sections of $\mc{IS^*}$ need not have compact support,  while the images of sections from $KS^*$ must.

\begin{lemma}\label{L: global mono}
$IS_0^*(X)=KS_0^*(X)=0$, i.e. there are no non-zero global presheaf sections with empty support. 
\end{lemma}
\begin{proof}

Let $\xi\neq 0\in IS_0^*(X)=IC^{\infty}_{n-*}(X;\mc{G}_0)$, and suppose that $|\xi|$ is empty. This means that for each point $x\in X$, the image of $\xi$ in $\lim_{x\in U} IC^{\infty}_{n-*}(X, X-\bar U;\mc G_0)=0$. So for all $x\in X$, there is a neighbohood $U_x$ of $x$ such that the support of $\xi$  lies in $X-\bar U_x$. In particular, then, $|\xi|\subset \cap_{x\in X} (X-\bar U_x)=\emptyset$, contradicting the non-triviality of $\xi$. The same arguments hold for $KS^*_0(X)$.
\end{proof}

\begin{lemma}\label{L: conjunctive}
$IS^*$  is conjunctive for coverings.
\end{lemma}
\begin{proof}
This proof is essentially the same as that in Swan for the sheaf of ordinary singular chains \cite[p. 85]{SW}. Let $X=\cup U_a$, and suppose $s_a\in IS^{n-i}(U_a)=IC^{\infty}_i(X, X-\bar U_a;\mc G_0)$ are such that $s_a|U_a\cap U_b= s_b|U_a\cap U_b$ for all $a,b$. We denote this common restriction by $s_{a\cap b}$. We need to find an $s\in IS^{n-i}(X)=IC^{\infty}_i(X;\mc G_0)$ such that $s|U_a=s_a$ for all $a$. 

Using the language of \cite[p. 85]{SW}, we will call a singular simplex $\sigma$  essential in an open set $U$ if $|\sigma|$  has a non-empty intersection with $\bar U$. We note once again (see Section \ref{S: IC}) that the coefficient of a singular simplex is determined entirely by the coefficient at  any interior point of $\Delta$ by the unique lifting property of local-coefficient systems, since $\sigma^{-1}(\Sigma)$ lies in $\bd \Delta$ so that $\Delta-\sigma^{-1}(\Sigma)$ is  contractible. 

Now, we claim that for any allowable singular $i$-simplex $\sigma:\Delta^i\to X$, the coefficient of $\sigma$ in $s_a$ is the same for all $a$ for which $\sigma$ is essential in $U_a$. If $\sigma$ is  essential in $U_a\cap U_b$, then the coefficients of $\sigma$ in $s_a$ and $s_b$ must agree since each restriction map $s_a\to s_{a\cap b}$ and $s_b\to s_{a\cap b}$ must preserve coefficients of simplices whose supports do not lie in $X-\overline{U_a\cap U_B}$.  If $\sigma$ is not essential in $U_a\cap U_b$, then since the $\{U_c\}$ form a covering and $|\sigma|$ is connected, there must be a finite  ``chain'' of elements of the covering $U_a=U_{c_0}, U_{c_1}, \ldots, U_{c_m}=U_b$ such that $U_{c_k}\cap U_{c_{k+1}}$ is non-empty for each $k$ and $\sigma$ is essential in each $U_{c_k}$ and each $U_{c_k}\cap U_{c_{k+1}}$. Inductively, the coefficients of $\sigma$ agree in all $s_{c_j}$ and $s_{c_j\cap c_{j+1}}$, so they agree in $s_a$ and $s_b$. 

We  then define $s= \sum_j g_j \sigma_j$, the sum over all singular simplices, where $g_j$ is the coefficient of $\sigma_j$ in $s_a$ for any $U_a$ in which $\sigma_j$ is essential. The arguments of the previous paragraph show that $g_j$ is well-defined. To see that $s$ is locally-finite, note that for any $x\in X$, $x\in U_a$ for some $a$, and  only those $\sigma$ which are essential in $U_a$ can have support that intersects $U_a$. But all such $\sigma$ must be in $s_a$, which is locally-finite. So $s$ is locally-finite in a neighborhood of every point; hence it is locally-finite. 

To see that $s$ is an allowable chain, we first note that each $i$ simplex in $s$ must be allowable, since an $i$-simplex can have non-zero coefficient in $s$ only if it has non-zero coefficient in some $s_a$, and each $s_a$ is an allowable chain. It remains to show that $\bd s$ is allowable. Each $i-1$ simplex  in $\bd s$ is an $i-1$ face of some $i$-simplex $\sigma$ in $s$ with non-zero coefficient. Suppose that $\tau$ is essential in $U_a$ and hence that $\sigma$ is as well. If $\tau$ is allowable, there is no issue. If $\tau$ is not allowable, then the coefficient of $\tau$ in $\bd s_a$ must be $0$. But any other singular $i$-simplex in $s_a$ that has $\tau$ in its boundary is also essential in $U_a$, and each occurs in $s$ with the same coefficent as it does in $s_a$. Thus  since the coefficient of $\tau$ in $\bd s_a$ must be $0$,  the coefficient of $\tau$ in $\bd s$ must be $0$ as well. 
\end{proof}

We can now show that the intersection homology groups defined in Section \ref{S: chains} can be recovered from the sheaf complex $\mc{IS}^*$ provided $X$ is paracompact and of finite cohomological dimension with respect to the ring $R$ such that $\mc G$ is a system of $R$ modules.

\begin{corollary}\label{C: section homology}
Let $X$ be a paracompact Hausdorff filtered space of finite cohomological dimension with respect to the ring $R$. Let $\mc{G}$ be a local coefficient system of $R$-modules on $X-X^{n-1}$. Then $ IH^{\infty}_{n-*}(X;\mc {G}_0)\cong H^*(\Gamma(X; \mc{IS^*}))$.
\end{corollary}
\begin{proof}
By definition, $IH^{\infty}_{n-*}(X;\mc{G}_0)=H_*(IC_{n-*}(X;\mc{G}_0))=H^*(IS^*(X))$. Since $IS^*_0(X)=0$ by Lemma \ref{L: global mono} and  $IS^*$ is conjunctive for coverings by Lemma \ref{L: conjunctive},  $IS^* (X)\cong \Gamma(X; \mc{IS^*})$ by \cite[I.6.2]{Br}, since $X$ is paracompact. So $IH^{\infty}_{n-*}(X;\mc G_0)\cong H^* (\Gamma(X; \mc{IS^*}))$.
\end{proof}

We next show that $\mc{IS^*}$ is homotopically fine.

\begin{proposition}\label{P: hom fine}
Let $X$ be a Hausdorff filtered space with coefficients $\mc G_0$. The sheaf $\mc{IS^*}$ is homotopically fine.
\end{proposition}
\begin{proof}
Let $\mc{U}=\{U_k\}$ be a locally-finite covering of $X$. We may impose  a well-ordering on $\mc U$. We must show that there exist endomorphisms $1_k$ and $D$ of $\mc{IS^*}$ such that $|1_k|\subset \bar U_k$ and $\sum 1_k=\text{id} -\bd D -D\bd$, i.e. $\sum 1_k$ is chain homotopic to the identity. The $1_k$ need not be chain maps. 

We first define a map $f_k: IC_*^c(X;\mc G_0)\to IC_*^c(U_k;\mc {G}_0)$ as follows: if  $\xi \in IC^c_i(X;\mc G_0)$, let $f_k(\xi)=\xi_k$ as defined in the proof of Proposition \ref{P: U small}. In other words, applying the subdivision operator $T$ of Proposition \ref{P: U small}, we take $T(\xi)$ and then discard from the chain $T(\xi)$ all $i$-simplices $\sigma$ for which $\psi(\sigma)\neq k$, where $\psi$ is also as defined in Proposition \ref{P: U small}. 
We must show that $f_k$ is a well-defined homomorphism of intersection chains. The image of each chain  $\xi \in IC_*^c(X;\mc G_0)$ under $f_k$  is an allowable intersection chain with support in $U_k$ by construction. $f_k$ is a homomorphism since it is determined linearly from what it does on singular simplices. 
Furthermore, if $j_k:IC^c_*(U_k;\mc {G}_0)\into IC^c_*(X;\mc {G}_0)$ is the inclusion  and $g_k=j_kf_k$, then $ \sum g_k=\iota T:IC^c_*(X;\mc {G}_0)\to IC^c_*(X;\mc {G}_0)$ is chain homotopic to the identity by a chain homotopy $D$ by Proposition \ref{P: U small}. 

Now, each map $g_k$ induces an  endomorphism of the presheaf $KS^*$ since each $KS^*(V)$ is a quotient of $IC_{n-*}(X;\mc G_0)$. On passing to sheaves, $g_k$ induces a map of sheaves $1_k$. Furthermore, the chain homotopy $D$ also descends to a map $\mc{D}$ of sheaves.

Let us check that these $1_k$ satisfy the desired properties; we  generalize the arguments of Swan \cite{SW} for ordinary singular homology. Suppose $x\notin U_k$, and, for any $s\in \mc{IS}^*_x$, let $\xi$ be a chain  in $IC_{n-*}^c(X, X-\bar V;\mc {G}_0)$ representing $s$, where $V$ is a neighborhood  of $x$. 
Since $|g_k(\xi)|$ is a compact subset of $ U_k$, we can find another neighborhood $W$ of $x$ such that $W\subset V$ and $\bar W\cap |g_k(\xi)|=\emptyset$.    Then $g_k(\xi)=0$ in $IC^c_*(X,X-\bar W)$. It follows that the image of $1_k$ is zero in the stalk $\mc{IS}^*_x$. Thus the support of $1_k$ must be in $U_k\subset \bar U_k$. 
As for the property that $\sum 1_k$ is homotopic to the identity, this follows from the fact that $\sum g_k$ is homotopic to the identity on presheaves; note that $\sum g_k$ is well-defined at the sheaf level since locally all but a finite number of terms are $0$. 
\end{proof}

Suppose $X$ is a paracompact Hausdorff filtered space. 
Since $\mc{IS^*}$ is homotopically fine by Proposition \ref{P: hom fine},    $H^*(H^p(X; \mc{IS^*}))=0$ for all $p>0$ by \cite[p. 172]{Br}.  So, as $\dim_R X<\infty$, there exists a spectral sequence with $E_2^{p,q}=H^p(X;\mc{H}^q(\mc{IS^*}))$ abuting to $H^{p+q}(\Gamma(X;\mc{IS^*}))=IH^{\infty}_{n-p-q}(X;\mc G_0)$ by \cite[IV.2.1]{Br} and Corollary \ref{C: section homology} (here $\mc{H}^*$ denotes  the derived cohomology sheaf). This is really just the hypercohomology spectral sequence for $\mc{IS}^*$, although we have to be a bit careful with our language as the term ``hypercohomology'' is often applied only to bounded below sheaf complexes, or at least those with bounded acyclic resolutions. However, our sheaf complex $\mc{IS^*}$ is not bounded below and a priori the acyclic ``resolution'' $L^*$ given in \cite[\S IV.1]{Br} (essentially the Cartan-Eilenberg flabby resolution) also will not be bounded below, though we will have $H^*(\Gamma(X;\mc{IS^*}))\cong H^*(\Gamma(X; L^*))$.  We will see, however, that if $X$ is a topological stratified pseudomanifold then we can in fact find a bounded below injective resolution $\mc{I^*}$ of $\mc{IS^*}$. It follows that $H^*(\Gamma(X;\mc{IS}^*))\cong H^*(\Gamma(X;\mc I^*)$   so that we can legitimately call $H^*(\Gamma(X;\mc{IS^*}))=IH_{n-*}^{\infty}(X;\mc G_0)$ the hypercohomology $\H^*(X;\mc{IS^*})$ of $\mc{IS^*}$ by any definition. We shall also continue to refer to $H^*(\Gamma(X;\mc{IS^*}))$ as the hypercohomology $\H^*(X;\mc{IS^*})$ as the spectral sequence does converge under our assumptions. 

\begin{remark}\label{R: extension} These observations provide a convincing argument that $\mc IS^*$ (or any of its representatives in the derived category) provides the correct ``sheafification'' of the intersection chain complex for any paracompact Hausdorff filtered space of finite cohomological dimension: the hypercohomology of this sheaf provides the intersection homology modules defined in Section \ref{S: chains}, and we will see below that on topological stratified pseudomanifolds, these modules agree with those defined by Goresky and MacPherson. Furthermore, on manifold weakly stratified spaces and for traditional perversities, the hypercohomology with compact supports gives the singular chain intersection homology studied by Quinn \cite{Q2}, who showed that the constant coefficient compactly supported singular intersection homology on such spaces is a topological invariant. 

What we lose on spaces more general than pseudomanifolds  is the  axiomatic characterization of the Deligne sheaf. Perhaps the  Goresky-MacPherson axioms can be extended to give an axiomatic description of $\mc IS^*$ on manifold weakly stratified spaces using the fact that points in such spaces have distinguished neighborhoods up to local stratum-preserving homotopy equivalence (see \cite{Q2}). However, it is only compact intersection homology theory that is a stratum-preserving homotopy invariant, and while locally-finite singular chain intersection homology may be a proper stratum-preserving homotopy type invariant, it is not evident that the  standard local stratum-preserving homotopy equivalences to distinguished neighborhood can be made proper. It is thus more difficult to compute closed support intersection homology of neighborhoods in such spaces, and it is unclear how to proceed with an analogue of the axiomatization as it is usually done on pseudomanifolds. Of course one could start with a coefficient system on $X-X^{n-1}$ and perform the Deligne construction, but it is not apparent that the hypercohomology of the resulting complex of sheaves will agree with any goemetric intersection homology  theory.
\end{remark}

\subsection{Restrictions to subspaces}\label{S: restriction}

In what follows we will also need to compare intersection homology sheaves $\mc{IS^*}$ on different spaces. Rather than use the notation $\mc{IS^*}(X)$ to indicate the space (which runs the danger of being confused with taking sections), we will indicate the space in subscript: $\mc{I_X S}^*$ being the sheaf of intersection chains on $X$. We continue to let $X$ be a paracompact Hausdorff space of stratified dimension $n$ and  to let the perversity or superperversity $\bar p$ and the coefficient system $\mc G_0$ remain fixed but absent from the notation. 

Suppose $W$ is an open subset of $X$, inheriting both the restricted filtration and coefficient system. We must study the relationship between the intersection chain sheaf on $W$ and the restriction to $W$ of the intersection chain sheaf on $X$.  In fact, they are quasi-isomorphic. 


\begin{proposition}\label{P: restriction}
Let $i: W\into X$ be an inclusion of an open subspace. Let $\mc{I_XS}^*$ and $\mc{I_WS^*}$  be the singular intersection chain sheaves on $X$ and $W$ respectively. Then there is a quasi-isomorphism $i^*\mc{I_XS^*}\to \mc{I_WS^*}$. 
\end{proposition}
\begin{proof}
We need to develop a map on chains that will induce  the quasi-isomorphism. We begin with a map of intersection chains  $r: IC^{\infty}_*(X)\to IC^{\infty}_*(W)$. We will define this map inductively over the dimension $j$ of the chains. As in the proof of Proposition \ref{P: U small}, we first define $r$ on  simplices of $C^{\infty}_*(X)$ and then indicate how to obtain a well-defined map of intersection chains.

For $j=0$, we define $r$ to be the restriction map that takes a singular $0$-simplex to itself if its support is in $W$ and to $0$ otherwise. Clearly this takes allowable $0$-chains to allowable $0$ chains. 

To define  $r$ on $1$-chains, we consider each singular $1$-simplex    $\sigma$ in $X$ and send it to a subdivision in $W$. By this we mean the following: If $\sigma: \Delta^1\to X$ is the given $1$-simplex, consider $\sigma^{-1}(W)$. This is an open subset of $\Delta^1$, which is a PL space, so we can find  a locally-finite PL-triangulation of $\sigma^{-1}(W)$, which we then replace with a singular triangulation based upon a partial ordering of vertices that respects the ordering of any vertices of $\Delta^1$ in $\sigma^{-1}(W)$. Our subdivision $r\sigma$ of $\sigma$ in $W$ is the chain consisting of the composition of $\sigma$ (as a chain map) with the singular simplices in this subdivision of $\sigma^{-1}(W)$. 
To see that $r\sigma$  is  locally-finite, let $x\in W$ and $V$ a neighborhood of $x$ with $\bar V\subset W$. Then $\overline{\sigma^{-1}(V)}\subset \sigma^{-1}(\bar V)\subset \sigma^{-1}(W)$. But $\sigma^{-1}(\bar V)$ is compact in $\Delta$ and so intersects only a finite number of simplices in the triangulation of $\sigma^{-1}(W)$. Thus only a finite number of singular simplices of $r\sigma$ have support intersecting $V$.  We define $r$ on $C_1(X)$ linearly by some choice of such triangulation on each basis $1$-simplex.  $r$ is clearly a chain map up to this point. $r$ then determines a map of intersection chains (also denoted $r$)  as in Proposition \ref{P: U small} by applying $r$ linearly to constituent simplices and subdividing coefficients in the obvious way.   The $1$-simplices in the image of an allowable chain under $r$ are allowable as in the proof of Lemma \ref{L: sub a}, and if $\xi$ is an intersection chain, $\bd (r\xi)=r(\bd \xi)$ is also allowable .

We now proceed inductively: assume $r$ defined on all $k$ chains for $k\leq j-1$, and let $\sigma:\Delta^j \to X$ be a singular $j$-simplex. This time we choose a locally-finite singular triangulation of $\sigma^{-1} (W)$ that agrees with the singular triangulation of $\sigma^{-1}(W)\cap \bd \Delta^j$ as determined by the induction hypothesis, as $\sigma|\bd \Delta$ determines a singular $j-1$ chain on which $r$ has already been defined. In particular, we triangulate $\sigma^{-1}(W)$ by singular chains using a partial ordering of vertices of a polyhedral subdivision consistent with that given by the subdivision of the boundary. $r$ is then extended linearly to all $j$-chains by choosing such a subdivision of each singular $j$-simplex. Again it is clear that $r$ will be a chain map, and  it induces a chain map $IC^{\infty}_*(X)\to IC^{\infty}_*(W)$ as in Proposition \ref{P: U small}.

Note that we are free to choose $r$ so that if $|\sigma|\subset W$, then  $r\sigma=\sigma$.

Now, if $U$ is an open set of $W$ with $\bar U\subset W$, the chain map $r$ induces a chain map $ IC_*^{\infty}(X,X-\bar U)\to IC_*^{\infty}(W, W-\bar U)$, and, in particular, if $x\in W$, we obtain maps $\lim_{x\in V} IC_*^{\infty}(X,X-\bar V)\to \lim_{x\in V} IC_*^{\infty}(W, W-\bar V)$. This induces a map from $i^*\mc{I_XS^*}=\mc{I_X S^*}|_W\to \mc {I_WS^*}$. To show that it is a quasi-isomorphism, we need only demonstrate an isomorphism on stalk cohomology. 

We first prove surjectivity. Let $s\in \mc H^*(\mc {I_WS^*})_x=H^*(\mc {I_WS}^*_x)$, which is isomorphic to both $ \lim_{x\in V} IH^c_{n-*}(W, W-\bar V)$ and $ \lim_{x\in V} IH^{\infty}_{n-*}(W, W-\bar V)$ since $\mc{I_WS^*}$ is the sheafification of both $I_WS^*$ and $K_WS^*$. 
Choose $V$ and  $[\xi]\in  IH^c_{n-*}(W, W-\bar V)$ such that $\xi$ represents $s$. Then 
 $[\xi]$ can be represented by a finite chain $\xi$ with $|\xi|\subset W$. We can also consider $\xi$ as an element of any of $IC^c_{n-*}(W)$, $IC^{\infty}_{n-*}(W)$, $IC^c_{n-*}(X)$, or $IC^{\infty}_{n-*}(X)$. Furthermore, since the inclusion $IH^c_{n-*}(W, W-\bar V)\to IH^{\infty}_{n-*}(W, W-\bar V)$ commutes with projection to the direct limit (which is essentially the content of the surjectivity half of Lemma \ref{L: same sheaf}), $s$ is also represented by $[\xi]\in IH^{\infty}_{n-*}(W, W-\bar V)$. Since $|\xi|\subset W$, $r\xi=\xi$, and it follows that the class of $\xi$ in $IH^{\infty}_{n-*}(X,X-\bar V)$ will map to $[\xi]\in IH^{\infty}_{n-*}(W, W-\bar V)$ under $r_*$. This shows that the sheaf map induced by $r$ is surjective on stalk cohomology. 

Next suppose that $s\in \mc H^*(i^*\mc {I_XS^*})_x\cong \lim_{x\in V} IH^c_{n-*}(X, X-\bar V)$ and that  $rs=0\in  H^*(\mc {I_WS^*})_x $. For each $\bar V\subset W$, it follows by excision (Lemma \ref{L: excision2}) that $IH^c_{n-*}(X, X-\bar V)\cong IH^c_{n-*}(W, W-\bar V)$ by excising the complement of $W$. Since the excision isomorphism is induced by inclusion, there is an open set $V$ and a finite chain $\xi$ in $IC^c_{n-*}(W)$ whose class in $IH^c_{n-*}(X, X-\bar V)$ represents $s$. As in the argument of the preceding paragraph, $\xi$ also represents $s$ as an element of $IH^{\infty}_{n-*}(X, X-\bar V)$. Furthermore, as a chain $r\xi=\xi$, and so $\xi$ also represents $rs \in \mc H^*(i^*\mc {I_WS^*})_x\cong \lim_{x\in V} IH^{\infty}_{n-*}(W, W-\bar V)\cong \lim_{x\in V} IH^c_{n-*}(W, W-\bar V)$. Now if $rs=0$, that implies that for some smaller open set $U$, $[\xi]=0\in IH^c_{n-*}(W, W-\bar U)$. Let $\Xi$ be a chain representing an element of $IC^c_{n-*}(W, W-\bar U)$ such that $\bd \Xi=\xi + \gamma$, where $|\gamma|\subset W-\bar U$. Now again, we can also think of $\Xi\in IC^{\infty}_{n-*}(X)$ and the equation  $\bd \Xi=\xi+\gamma$ continues to hold in this module. Thus $[\xi]=0\in IH^{\infty}_{n-*}(X, X-\bar U)$, and $s=0$. 
\end{proof}

We can also construct an explicit quasi-isomorphism in the other direction from $\mc{I_WS}^*$ to $\mc{I_XS}^*|_W$. In fact, this quasi-isomorphism is much easier to construct, though we will need the other one explicitly in what follows. 
Recall that these sheaves are the sheafifications of the pre-sheaves $K_WS^*$ and $K_XS^*$, where $K_XS^{n-*}(V)=IC^c_*(X, X-\bar V; \mc G_0)$, and similarly for $K_WS^{n-*}(V)$. Since $W$ is open in $X$, $\mc{I_XS}^*|_W$ is simply the sheafification of the restriction of $I_XS^*$ to $W$. But for $V$ open with $\bar V\subset W$, the  map induced by inclusion $i: IC^c_*(W, W-\bar V; \mc{G}_0)\to IC^c_*(X, X-\bar V; \mc{G}_0)$ induces a map of presheaves which, by Lemma \ref{L: excision2}, induces a homology isomorphism by excision of $X-W$. Since for all pairs $U\subset V$ the diagram
\begin{equation*}
\begin{CD}
IH^c_*(W, W-\bar V; \mc{G}_0)&@>i_*>>& IH^c_*(X, X-\bar V; \mc{G}_0)\\
@VVV&&@VVV\\
IH^c_*(W, W-\bar U; \mc{G}_0)&@>i_*>>& IH^c_*(X, X-\bar U; \mc{G}_0)
\end{CD}
\end{equation*}
commutes,
the map of direct systems induces an isomorphism $\lim_{x\in V}  IH^c_*(W, W-\bar V; \mc{G}_0)\to \lim_{x\in V}IH^c_*(X, X-\bar V; \mc{G}_0)$. This shows that $i$ induces  a quasi-isomorphism of sheaves. 

Furthermore, $r$ and $i$ are quasi-inverses: Let $x\in W$ and let $s\in \mc{I_WS}^*_x$. It follows as in the proof of Proposition \ref{P: restriction} that $s$ can be represented locally by a finite chain $\xi\in IC^{c}_*(W, W-\bar V)\subset IC^{\infty}_*(W, W-\bar V)$ for some neighborhood $V$ of $x$. Under $i$, $\xi$ also represents an element of $IC_*^{\infty}(X, X-\bar V)$, and $r(\xi)=\xi$ since $|\xi|\subset W$. Thus $ri(\xi)$ also represents $s$.
So $ri$ is the identity at each stalk in $W$, and  it induces the identity isomorphism on cohomology stalks. Now, since $r_*i_*$ is the identity on cohomology stalks, we have  $ i_*r_*i_*=i_*$. But $i_*$ is a homology isomorphism, so $i_*r_*$ is also the identity map on cohomology stalks. Thus $i$ and $r$ are quasi-inverses.

\smallskip
We next utilize Proposition \ref{P: restriction},  together with our earlier intersection homology computations,  to compute the map induced by restriction from  the intersection homology of a distinguished neighborhood to that of the corresponding deleted distinguished neighborhood. This will be important below in demonstrating that the sheaf attaching map is an isomorphism in a certain range.

\begin{proposition}\label{P: cone hom2}
Let $L^{k-1}$ be a compact filtered space.  Let $x$ be the cone point of $cL$.  The restriction map $r: IC_*^{\infty}( cL\times \R^{n-k})\to IC_*^{\infty}((cL-x)\times \R^{n-k})$ of Proposition \ref{P: restriction} induces an isomorphism on homology in dimensions $\geq n-\bar p(k)$ and the $0$ map otherwise. 
\end{proposition}
\begin{proof}
First note that it follows from Propositions \ref{P: dist ngbd} and \ref{P: cone} that $IH_*^{\infty}( cL\times \R^{n-k})=0$ for $*<n-\bar p(k)$.

By Propositions \ref{P: dist ngbd} and \ref{P: hole}, the homology modules  $IH_*^{\infty}( cL\times \R^{n-k})$ and $IH_*^{\infty}((cL-x)\times \R^{n-k})$ are respectively isomorphic to $IH_{*-(n-k)}^{\infty}( cL)$ and $ IC_{*-(n-k)}^{\infty}(cL-x)$. By Proposition \ref{P: cone} the former is isomorphic to $IH_{*-(n-k)-1}(L)$ for $*-(n-k)\geq k-\bar p(k)$ (i.e. for $*\geq n-\bar p(k)$), while by Proposition \ref{P: times R}, the latter is always isomorphic to $IH_{*-(n-k)-1}(L)$. Thus abstractly the modules are isomorphic in the appropriate range. We must show that the isomorphism is induced by $r$. 

By the Remarks following the proofs of Propositions \ref{P: times R}, \ref{P: cone}, \ref{P: dist ngbd}, and \ref{P: hole}, if $\xi$ is a chain representing an element $[\xi]\in IH_{*-(n-k)-1}(L)$, then the image of $[\xi]$ in  $IH_*^{\infty}( cL\times \R^{n-k})$ and $ IH_*^{\infty}((cL-x)\times \R^{n-k})$ under these isomorphisms (in the appropriate dimension ranges) are given, respectively, by the chains $c\xi \times \R^{n-k}$ and $\xi \times (0,1)\times \R^{n-k}$ (where we have identified $cL-x$ with $L\times (0,1)$). 

So now let $[\zeta] \in IH_i^{\infty}( cL\times \R^{n-k})$,   $i\geq n-\bar p(k)$, and let $\xi$ represent the corresponding class in $IH_{i-(n-k)-1}(L)$ so that we can take $c\xi\times \R^{n-k}$ as a chain representing the class $[\zeta]$. It will suffice to show that $r(c\xi\times \R^{n-k})$ represents the correct corresponding class in $IH_i^{\infty}((cL-x)\times \R^{n-k})$, i.e. that it can be written as $\xi\times (0,1)\times \R^{n-k+1}$. 
But this follows now from the definitions:

If $\sigma_j:\Delta_j\to L$ is a singular simplex in $\xi$, each  $c\sigma_j\times \R^{n-k}$ of $c\xi\times \R^{n-k}$ comes from composing a singular triangulation of $c\Delta_j \times \R^{n-k}$ with the product of $\sigma_j$ and the identity maps in the $\R^{n-k}$ and cone directions.
In other words, if $(x,t,s)\in c\Delta_j\times \R^{n-k}$, with $x\in \Delta_j$, $t\in (0,1)$, and $s\in \R^{n-k}$, then, as a map, $(c\sigma_j\times \R^{n-k})(x,t,s)=(\sigma(x),t,s)\in cL\times \R^{n-k}$. Composing with the singular triangulation of $c\Delta_j\times \R^{n-k}$ gives 
$c\sigma_j\times \R^{n-k}$ as a chain. Then of course $c\xi\times \R^{n-k}$ is the sum of the $c\sigma_j\times \R^{n-k}$ weighted by the similarly treated coefficients $cn_j\times \R^{n-k}$. 

Since $|c\sigma_j\times \R^{n-k}|$ intersects $x\times \R^{n-k}$ only along the image of $z\times \R^{n-k}$, where $z$ is the cone point  of $c\Delta_j$, the effect of $r$ on $c\sigma_j\times \R^{n-k}$ is, by definition, the chain obtained  by retriangulating $(c\Delta_j-z) \times \R^{n-k}$ according to the construction of $r$ in Proposition \ref{P: restriction} and composing with the map $c\sigma_j\times \R^{n-k}$ restricted to $(c\Delta_j-z) \times \R^{n-k}$. 
But of course, since $r$ is a chain map, these give subdivisions compatibile among the  $\Delta_j\times (0,1)\times \R^{n-k}$, which then can be used to define $\xi \times (0,1)\times \R^{n-k}$. But we know that $\xi\times (0,1)\times \R^{n-k}$ represents the desired class in  $IH_*^{\infty}((cL-x)\times \R^{n-k})$.
\end{proof}

\subsection{Agreement of singular $IH$ theory with Deligne-sheaf $IH$ theory}\label{S: axioms}

In this section, we show that on a paracompact topological pseudomanifold, the sheaf complex $\mc{IS}^*$ defined in Section \ref{S: sheafify} is quasi-isomorphic to the corresponding Deligne sheaf. This implies that, for traditional perversities, $IH^{\infty}_*(X;\mc G_0)$ as defined in Section \ref{S: chains} is isomorphic to the intersection homology of Goresky and MacPherson \cite{GM2}, while for superperversities with $\bar p=1$, it is isomorphic to the superperverse intersection homology occuring in the superduality theorem of Cappell and Shaneson \cite{CS}. 

Let $X$ be a paracompact $n$-dimensional topological stratified pseudmanifold. Then, as noted in \cite[p. 60]{Bo}, $X$ has cohomological dimension $n$ and is locally compact, hence locally paracompact. In addition, distinguished neighborhoods are each paracompact: by general topology (see \cite[Th. 2-65]{HY}),  it suffices that each distinguished neighborhood $N$ is locally compact Hausdorff and the union of a countable number of compact sets, which is easily verified. 
We  also  fix a ground ring $R$ with unit and of finite cohomological dimension.  Let $\mc {G}$ be a local coefficient system of $R$ modules on $X-\Sigma=X-X^{n-2}$.

Let $\mc {P^*}$ denote the Deligne sheaf on $X$ determined by a traditional perversity or superperversity $\bar p$ and the local coefficient system $\mc {G}$ on $X-\Sigma$. When we wish to emphasize the coefficient system, we will write $\mc P^*_{\mc G}$. We recall that this is the sheaf on $X$ defined inductively as follows (see \cite{GM2} or \cite[\S V.2]{Bo}): On $X-X^{n-2}$, $\mc{P}^*_2=\mc{G}$. Then for each $k\geq 2$, $\mc{P^*}_{k+1}=\tau_{\leq p(k)}Ri_{k*}\mc{P}^*_k$, where $i_k:X-X^{n-k}\to X-X^{n-k-1}$ is the open inclusion, $Ri_{k*}$ is its right derived functor, and $\tau_{\leq p(k)}$ is the truncation. Then $\mc P^*=\mc P^*_{n+1}$, and  the  intersection homology of $X$ as defined in \cite{GM2} is isomorphic to the hypercohomology $\H^*(\mc P^*)$. 

Our main theorem essentially says that our complex of intersection chains with coefficients $\mc G_0$ gives the same intersection homology modules as those obtained from the hypercohomology of $\mc {P^*}$, for either perversities or superperversities. More specifically, we will see that  $I^{\bar p}H^{\infty}_{n-*}(X;\mc G_0)\cong \H^*(\mc P^*_{\mc G\otimes_R \mc O})$, where $\mc O$ is the orientation $R$-module on the manifold $X-\Sigma$. So if $X-\Sigma$ is orientable, $I^{\bar p}H^{\infty}_{n-*}(X;\mc G_0)\cong \H^*(\mc P^*_{\mc G})$, and even if it is not orientable,  we have $\H^*(\mc P^*_{\mc G})\cong I^{\bar p}H^{\infty}_{n-*}(X;(\mc G\otimes_R \mc O)_0)$, since $\mc G \otimes_R \mc O \otimes_R \mc O\cong \mc G$.  

\begin{theorem}\label{T: main theorem}
Let $X$ be a paracompact $n$-dimensional topological stratified pseudomanifold, and let $\mc G$ denote a local coefficient system of $R$ modules on $X-\Sigma$. Then $IH^{\infty}_{n-i}(X;\mc G_0)\cong \H^i(\mc {P}^*_{\mc G\otimes_R \mc O})$. In particular, if $X-\Sigma$ is orientable, then $IH^{\infty}_{n-i}(X;\mc G_0)\cong \H^i(\mc {P}^*_{\mc G})$, and if $X$ is compact, then $\H^i(\mc {P}^*_{\mc G\otimes_R \mc O})\cong IH^{c}_{n-i}(X;\mc G_0)$.
\end{theorem}

\begin{proof}

We will proceed by induction on the dimension $n$ of the pseudomanifold $X$. 

If $X$ has dimension $0$, then $X$ is a collection of discrete points $X=\amalg x_j$, and $\mc{G}$ is a collection of modules $G_j$, one for each point. The singular intersection homology reduces to ordinary locally-finite singular homology and so $IH_*^{\infty}(X;\mc G_0)=H_{*}^{\infty}(X; \mc G)$, which is $\prod G_j$ in dimension $*=0$ and $0$ otherwise. This of course agrees with $\H^{*}(\mc {P_{\mc G}})$, which is $\Gamma(X,\mc G)=\prod G_j$ in dimension $0$ and $0$ otherwise. 

So now, inductively, we assume the theorem has been proven for all stratified pseudomanifolds of dimension $<n$, and we fix for the remainder of the argument a perversity or superperversity $\bar p$, a stratified topological pseudomanifold $X$ of dimension $n$, and a coefficient system $\mc G$ on $X-\Sigma$. All sheaves from here out will be with respect to these fixed choices, omitted from the notation.

We first truncate to a quasi-isomorphic sheaf complex that has the benefit of being bounded.

\begin{lemma}\label{L: 0 dims}
The truncation $\mc{IS^*}\to \tau^{\geq c} \mc{IS^*}$ is a quasi-isomorphism for any $c\leq 0$. 
\end{lemma}
\begin{proof}
Here $\tau^{\geq c}$ is the standard truncation functor that gives
\begin{equation*}
(\tau^{\geq c}\mc{IS})^i=
\begin{cases}
0, & i<c\\
\cok(d_{i-1}), & i=c\\
\mc{IS}^i, & i>c.
\end{cases}
\end{equation*}
Truncation always induces a quasi-isomorphism for $i\geq c$. We must show that the derived cohomology sheaf $\mc{H}^i(\mc{IS^*})=0$ for $i<0$. 

$\mc{H}^i(\mc{IS^*})$ is the sheafification of the presheaf $U\to H^i(IS^*(U))$, which is $U\to IH_{n-i}^{\infty}(X,X-\bar U;\mc {G}_0)$. To compute stalks, we take the direct limit over distinguished neighborhoods of the point $x$. By Corollary \ref{C: essentially constant}, which says that the direct system is essentially constant, it suffices to fix a distinguished neighborhood of the form $N_{\alpha}$, $0<\alpha<1$, (in the notation of Lemma \ref{L: rest dist}) and to show that $IH_{n-i}^{\infty}(X,X-\bar N_{\alpha};\mc {G}_0)=0$ for $i<0$. Let $L$ denote the compact $n-j-1$ dimensional link pseudmanifold of $x$ so that $N_{\alpha}\cong cL\times \R^j$.

Now  by Lemma \ref{L: cocompactness}, $IH^{\infty}_*(X,X-\bar N_{\alpha};\mc G_0)\cong IH^c_*(X,X-\bar N_{\alpha};\mc G_0)$, which by stratum-preserving homotopy equivalence is isomorphic to $IH^c_*(X,X-x;\mc G_0)$. By excision (Lemma \ref{L: excision2}), this is isomorphic to $IH^c_*(N_{\alpha},N_{\alpha}-x;\mc G_0)$, which again by  Lemma \ref{L: cocompactness} is isomorphic to $IH^{\infty}_*(N_{\alpha},N_{\alpha}-x;\mc {G}_0)\cong IH^{\infty}_*(cL\times \R^j,cL\times \R^j-x;\mc {G}_0)$. It was seen in the proof of Proposition \ref{P: dist ngbd} that $IH^{\infty}_*(cL\times \R^j,cL\times \R^j-x;\mc {G}_0)\cong IH^{\infty}_*(cL\times \R^j;\mc {G}_0)$, and 
by Propositions \ref{P: dist ngbd} and \ref{P: cone}, $IH^{\infty}_*(cL\times \R^j;c\mc G_0\times \R^j)\cong IH^{\infty}_{*-j}(cL;c\mc G_0)\cong  IH^{\infty}_{*-j-1}(L;\mc G_0)$ if $*-j\geq n-j-\bar p(n-j)$ and $0$ otherwise.

Now $L$ has dimension $n-j-1<n$, so by the induction hypothesis $IH_*(L;\mc G_0)\cong \H^{n-j-1-*}(\mc P^*_{L,\mc G\otimes_R \mc O})$,  where $\mc P^*_{L,\mc G\otimes_R \mc O}$ is the Deligne sheaf on $L$ induced by the restricted coefficient system and orientation sheaf  on $L$. In particular, $IH_*(L;\mc G_0)=0$ in dimensions  $> n-j-1$. So $IH^{\infty}_*(X, X-\bar N_{\alpha};\mc G_0)=0$ for $*>n$. Thus $IH^{\infty}_{n-*}(X, X-\bar N_{\alpha};\mc G_0)=0$ for $*<0$,and so $\mc{H}^i(\mc{IS^*})=0$ for $i<0$. 
\end{proof}

So we have a quasi-isomorphism $\mc {IS^*}\to \tau^{\geq 0} \mc{IS^*}$. As a bounded from below complex, $\tau^{\geq 0} \mc{IS^*}$ has an injective resolution, say $\tau^{\geq 0}\mc{IS^*}\to \mc I^*$, such that $I^j=0$ for $j<0$.  Then $H^*(\Gamma(I^*))$ would be the standard hypercohomology  $ \H^*(\tau^{\geq 0} \mc{IS^*})$. But now the composition  $\mc {IS^*}\to \tau^{\geq 0} \mc{IS^*}\to \mc I^*$ of quasi-isomorphisms is a quasi-isomorphism. Therefore, since $\mc{IS^*}$ is homotopically fine, $\mc I^*$ is injective, and $\dim_R(X)<\infty$, the associated map of spectral sequences gives an isomorphism  $H^i(\Gamma(\mc {IS^*}))\to H^i (\Gamma(I^*))=\H^i(I^*)$ for all $i$ (see \cite[IV.2.2]{Br}). But we already know by Lemma \ref{C: section homology} that $H^i(\Gamma(\mc {IS^*}))=H_{n-i}(X;\mc G_0)$. This argument legitimitizes our earlier claim  that $H^i(\Gamma(\mc {IS^*}))$ should be called the hypercohomology of $\mc{IS^*}$. 
It now suffices to show that  $\tau^{\geq 0} \mc{IS}^*$ (and hence $\mc I^*$) is quasi-isomorphic to the Deligne sheaf and hence yields the same hypercohomology.  

For convenience of notation, we let $\mc T^*=\tau^{\geq 0} \mc{IS^*}$.

To verify that a complex of sheaves is quasi-isomorphic to $\mc{P}_{\mc G\otimes_R \mc O}$, it is only necessary to check the intersection sheaf axioms  (see \cite{GM2}, \cite[\S V.2]{Bo}). In particular, we must show that 
\begin{enumerate}
\item $\mc T^*$ is bounded and  $\mc T^*=0$ for $*<0$,
\item $\mc T^*|_{X-\Sigma}$ is quasi-isomorphic to $\mc G\otimes_R \mc O$,
\item for $x\in X_{n-k}-X_{n-k-1}$, $\mc H^j(\mc T^*)_x=0$ if $j>\bar p(k)$
\item the attaching map $\alpha: \mc T^*|_{X-X_{n-k-1}}\to Ri_{k_*}\mc T^*|_{X-X_{n-k}}$ is a quasi-isomorphism up to dimension $\bar p(k)$, where $i_k:X-X_{n-k}\to X-X_{n-k-1}$ is the inclusion.
\end{enumerate}

We show that these axioms are satisfied:

1. $\mc T^*=\tau^{\geq 0} \mc{IS^*}$ is certainly bounded below due to the truncation. It is also bounded above since there are no singular chains of negative dimension.

2. Consider the restriction of $\mc{IS^*}$ to $U_2=X-\Sigma$. By Proposition \ref{P: restriction}, $\mc{IS^*}|_{U_2}$ is equal to the sheaf $\mc{I}_{U_2}\mc S^*$ of intersection chains on $U_2$. But $U_2$ is an $n$-manifold, so its  intersection chain sheaf agrees with the ordinary locally-finite chain sheaf $\mc{S}^*$ on $U_2$ with the given coefficients $\mc{G}$. Hence at $x\in U_2$, $\mc {H}^*(\mc{S^*})=\lim_{x\in V} H^{\infty}_{n-*}(X, X-V;\mc {G})$. Stalkwise, this is $0$ for $*\neq 0$, and in dimension $0$ we obtain the stalk $G$ of $\mc G$, the local orientation determining the identification. Globally, we obtain $\mc {G}\otimes_R \mc{O}$, the tensor product with the orientation sheaf.  So $\tau^{\leq 0}\mc{IS^*}|_{U_2}$ is quasi-isomorphic to $\mc{G}\otimes_R \mc O$.

3. We must show that for all $x\in S_{n-k}=X_{n-k}-X_{n-k-1}$,  $\mc{H}^j(\mc{T}^*)_x=0$ for $j>\bar p(k)$. So let $x\in S_{n-k}$. Then $x$ has a cofinal system of distinguished neighborhoods of the form $N\cong cL\times \R^{n-k}$, where $L$ is a compact $k-1$ pseudomanifold. So we look at $\lim_{x\in V} IH^{\infty}_{n-j}(X, X-\bar V;\mc {G}_0)\cong H^j(\mc{T}^*)_x$. By the arguments in Lemma \ref{L: 0 dims}, for sufficiently small $N$,
$I^{\bar p}H^{\infty}_*(X,X-\bar N;\mc G_0)\cong   IH^{\infty}_{*-(n-k)-1}(L;\mc G_0)$ if $*-(n-k)\geq k-\bar p(k)$ and $0$ otherwise. So this module is $0$ if $n-*>\bar p(k) $, which implies that $\mc{H}^j(\mc{T}^*)_x=0$ for $j>\bar p(k)$.

4. Lastly, we must show that the attaching map $\mc{T^*}|_{U_{k+1}}\to Ri_{k*} (\mc{T^*}|_{U_k})$ is a quasi-isomorphism up to dimension $\bar p(k)$, where $i_k$ is the inclusion $U_k=X-X_{n-k}\into X-X_{n-k-1}=U_{k+1}$. This map is automatically a stalk quasi-isomorphism  at points of $U_k$, so it is only necessary to check points in $U_{k+1}-U_k=S_{n-k}=X_{n-k}-X_{n-k-1}$. As noted in \cite[p.50]{Bo}, demonstrating this isomorphism amounts to showing that for a system of distinguished neighborhoods $N$ of $x\in S_{n-k}$,  $H^i(\mc{T}^*_x)\cong \lim_{x\in N} \H^i(N-N\cap S_{n-k};\mc{T}^*)$. 

Now $H^i(\mc{T}^*_x)\cong \mc{H}^i(\mc T)_x$, and we have already seen that this is equal to  $I^{\bar p}H^{\infty}_{(n-i)-(n-k)-1}(L;\mc G_0)$ if $n-i\geq n-\bar p(k)$ and $0$ otherwise. Meanwhile, $\H^i(V-V\cap S_{n-k};\mc{T}^*)$ is, by definition, the cohomology of $\Gamma(V-V\cap S_{n-k};\mc I^*)$, where $\mc I^*$ is an injective resolution of $\mc T^*$. Since restriction is an exact functor and the restriction of an injective sheaf to an open set is injective, $\mc I^*|_{V-V\cap S_{n-k}}$ is an injective resolution of $\mc{T^*}|_{V-V\cap S_{n-k}}$. But by Proposition \ref{P: restriction}, for an open set $W$, $\mc{IS^*}|_{W}$ is quasi-isomorphic to $\mc{I_WS^*}$, the intersection chain sheaf on $W$. Since the functor $\tau^{\geq 0}$ commutes with restriction, we see that $\mc{T^*}|_{V-V\cap S_{n-k}}$ is quasi-isomorphic to $\mc{T}^*_{V-V\cap S_{n-k}}$, and $\mc I^*|_{V-V\cap S_{n-k}}$ is an injective resolution of $\mc T^*_{V-V\cap S_{n-k}}$ (where we use $\mc T^*_Z$ to denote $\tau^{\geq 0}\mc{I_ZS^*}$; N.B. this is not the same as the restriction to $Z$, $\mc T^*|_Z$). So we have isomorphisms 
\begin{align*}
\H^i(V-V\cap S_{n-k};\mc{T}^*)& \cong H^i(\Gamma(V-V\cap S_{n-k};\mc I^*))\\
&\cong  H^i(\Gamma(V-V\cap S_{n-k};I^*|_{V-V\cap S_{n-k}}))\\
&\cong \H^i (V-V\cap S_{n-k} ;T^*_{V-V\cap S_{n-k}})\\
&\cong IH^{\infty}_{n-i}(V-V\cap S_{n-k};\mc G_0).
\end{align*} 
Now $V-V\cap S_{n-k}\cong  (cL-x)\times \R^{n-k}$, and by Propositions \ref{P: hole} and \ref{P: times R}, $IH^{\infty}_{n-i}((cL-x)\times \R^{n-k};c\mc G_0\times \R^k)\cong IH^{\infty}_{(n-i)-(n-k+1)}(L;\mc G_0)$. 

 So we see that, abstractly, $H^i(T^*_x)\cong \H^i(V-V\cap S_{n-k};\mc{T})$ if $n-i\geq n-\bar p(k)$, i.e. if $i\leq \bar p(k)$. 
We need to  show that this isomorphism is indeed induced by the attaching map. This will be done in the following proposition, which will complete our proof of the Theorem \ref{T: main theorem}. 
\end{proof}

\begin{proposition}
The attaching map $\alpha: \mc{T^*}|_{U_{k+1}}\to Ri_{k*} (\mc{T}|_{U_k})$ is a quasi-isomorphism up to dimension $\bar p(k)$. 
\end{proposition}
\begin{proof}

Since we are only concerned with quasi-isomorphisms, we can replace $\mc {T^*}|_{U_{k+1}}$ by its injective resolution $\mc{I^*}$:

\begin{diagram}
\mc{T^*}|_{U_{k+1}} & \rTo & i_{k*}i_k^* (\mc{T}|_{U_{k+1}})& \rTo &Ri_{k*}i_k^* (\mc{T}|_{U_{k+1}})\\
\dTo &&\dTo& \ruTo^{=}\\
\mc{I^*} &\rTo&i_{k*}i_k^* \mc{I^*}
\end{diagram}
The two vertical maps are induces by the injective resolution $\mc{T^*}|_{U_{k+1}}\to \mc{I^*}$. The top row is the row we wish to show 
is a quasi-isomorphism in the desired range. The box commutes by functoriality. The triangle commutes by definition. The lefthand map is a quasi-isomorphism since it is an injective resolution. So it suffices to show that the bottom map is a quasi-isomorphism.
This is automatic  at all points of $U_{k+1}-U_k=S_{n-k}$;  we must check stalk maps  at points  $x\in S_{n-k}$. 

Let $x\in S_{n-k}$. We want to show that
$\lim_{x\in V} \Gamma(V, \mc{I^*})\to \lim_{x\in V} \Gamma(V, i_{k*}i_k^*\mc{I^*})=\lim_{x\in V} \Gamma(V-V\cap S_{n-k}, i_k^*\mc{I^*})=\lim_{x\in V}\Gamma(V-V\cap S_{n-k}, \mc{I^*})$ induces  an isomorphism on cohomology in the appropriate range. We will show that, in fact, this map of sections $\Gamma(V,\mc I^*)\to \Gamma(V-V\cap S_{n-k}\mc I^*)$ (induced my restriction) induces cohomology isomorphisms for each distinguished neighborhood $V$. As the distinguished neighborhoods constitute a cofinal system of neighborhoods and restriction of sections is natural, this will induce the quasi-isomorphism on the direct limits. 

So let $V$ be a distinguished neighborhood of $x$, and let $Z=V\cap S_{n-k}$. 
Since 
\begin{equation*}
\begin{CD}
 \Gamma(V, \mc{I^*})&@>>>& \Gamma(V-Z, \mc{I^*})\\
@V=VV&&@V=VV\\
\Gamma(V, \mc{I^*}|_V)&@>>>& \Gamma(V-Z, \mc{I^*}|_V)
\end{CD}
\end{equation*}
commutes, where the horizontal arrows are restrictions, we can limit attention to $\mc{I^*}|_V$, which, by restriction, is an injective resolution of $\mc T^*|_V$ and $\mc{IS^*}|_V$. Since we have a quasi-isomorphism $i:\mc{I_VS}^*\to \mc{IS}^*|_V$ by Proposition \ref{P: restriction}, we can also consider $\mc{I^*}|_V$ as an injective resolution of $\mc{I_VS}^*$, and we obtain the following commutative diagram induced by these resolutions:
\begin{equation}\label{CD1}
\begin{CD}
\Gamma(V, \mc{I^*}|_V)&@>>>& \Gamma(V-Z, \mc{I^*}|_V)&@>=>>&\Gamma(V-Z,\mc{I^*}|_{V-Z})\\
@AAA&&@AAA &&@AAA\\
\Gamma(V, \mc{I_VS^*})&@>>>& \Gamma(V-Z, \mc{I_VS^*})&@>=>>&\Gamma(V-Z, \mc{I_VS^*}|_{V-Z})&.
\end{CD}
\end{equation}

Since $\mc{I_VS^*}$ is the sheafification of the presheaf $I_VS^*$ and since $V$ is paracompact, $\Gamma(V, \mc{I_VS^*})\cong IC^{\infty}_{n-*}(V)$ by 
Lemmas \ref{L: global mono} and  \ref{L: conjunctive} and by  \cite[I.6.2]{Br}.
Similarly, $\Gamma(V-Z, \mc{I_{V-Z}S^*})\cong IC^{\infty}_{n-*}(V-Z)$. 
We consider then the diagram
\begin{diagram}[LaTeXeqno]\label{CD2}
\Gamma(V, \mc{I_VS^*}) & \rTo^R & \Gamma(V-Z, \mc{I_VS^*}|_{V-Z})\\
 &&\dTo^{\bar r}\\
\uTo^{\cong}_{\phi}&& \Gamma(V-Z,\mc{I_{V-Z}S^*})\\
&&\uTo^{\cong}_{\psi}\\
 IC^{\infty}_{n-*}(V)&\rTo^r &IC^{\infty}_{n-*}(V-Z)&,
\end{diagram}
where the map $r$ is the restriction map of Proposition \ref{P: restriction},  $\bar r$ is the map it induces  on sheaves, $\phi$ and $\psi$ are the natural maps induced by sheafification, and $R$ is the composition of the bottom row of Diagram \eqref{CD1}.  To see that this diagram commutes, let us begin with a chain $\xi\in  IC^{\infty}_{n-*}(V)$. The image of $ \xi$ under $\phi$ is a section whose germ at $z\in V$ is the image of $\xi$ in $\lim_{z\in U} IC^{\infty}_{n-*}(V,V-\bar U)$. The image of $\phi(\xi)$ under $R$ is then simply the restriction of this section to the stalks in $V-X$. Finally, $\bar r$, since it is induced by the map $r$ on presheaves, takes $R\phi(\xi)$ to a section such that at each  $z\in V-Z$, the germ is $\lim_{z\in U} r\xi \in \lim_{z\in U} IC^{\infty}_{n-*}(V-Z,(V-Z)-\bar U)$. But this describes precisely the image under $\psi$ of $r\xi$. 

We now require one last commutative diagram:
\begin{diagram}[LaTeXeqno] \label{CD3}
\mc{IS^*}|_{V-Z} &\rTo & \mc{T}^*|_{V-Z} &\rTo & \mc{I^*}|_{V-Z}\\
\uTo^{\bar i} &&\uTo && \uTo\\
\mc{I_VS^*}|_{V-Z} &\rTo  &\mc{T_V^*}|_{V-Z}& \rTo & \mc{I_V^*}|_{V-Z}\\
\dTo^{\bar r} && \dTo && \dTo\\
\mc{I_{V-Z}S^*} & \rTo & \mc{T_{V-Z}^*}& \rTo & \mc{I_{V-Z}^*}
\end{diagram}
Here $\mc T^*_Y$ and $\mc I^*_Y$ are defined analogously to the sheaves $\mc T^*$ and $\mc I^*$, but over the space $Y$. These should not be confused with the restrictions, denotes $\mc T^*|_Y$ and $\mc I^*|_Y$. The map $\bar i$ is induced by the inclusion $i: K_VS^*\to KS^*$.
The left two columns commute by functoriality of the truncation functor. The maps in the right column exist and create commutative squares since the various $\mc{I}^*$ can each be taken as Cartan-Eilenberg resolutions, which can be completed to commutative diagrams \cite{CE} (note that the restriction of a Cartan-Eilenberg resolution to an open set is a Cartan-Eilenberg resolution of the restriction). The two lefthand maps are quasi-isomorphisms by Proposition \ref{P: restriction} and the discussion following it,
and the horizontal maps are  quasi-isomorphisms by Lemma \ref{L: 0 dims} and by construction.
Thus all maps in the diagram are quasi-isomorphisms. This diagram induces a commutative diagram on taking sections over $V-Z$. Putting the appropriate pieces of the commutative diagrams \eqref{CD1}, \eqref{CD2}, and \eqref{CD3} together, gives us a commutative diagram

\begin{diagram}
\Gamma(V, \mc{I^*}|_V)&\rTo &&& \Gamma(V-Z, \mc{I^*}|_V)=\Gamma(V-Z,\mc{I^*}|_{V-Z})\\
\uTo &&&\ruTo& \uTo\\
&&\Gamma(V-Z,\mc{I_VS^*}|_{V-Z})&& \Gamma(V-Z,\mc{I_V^*}|_{V-Z})\\
&\ruTo^{R}&\dTo_{\bar r}&&\dTo\\
\Gamma(V, \mc{I_VS^*}) && \Gamma(V-Z, \mc{I_{V-Z}S^*})&\rTo^f & \Gamma(V-Z, \mc{I_{V-Z}^*})\\
\uTo^{\phi}_{\cong}&&\uTo_{\psi}^{\cong} \\
IC^{\infty}_{n-*}(V) &\rTo^r &IC^{\infty}_{n-*}(V-Z)
\end{diagram} 
The upper right diagonal map is the composition of the upper left vertical map and the top row of Diagram \eqref{CD3}.
Now,  the two sides of the bottom square are quasi-isomorphisms, and the bottom map is a homology isomorphism for $*\leq \bar p(k)$ by Proposition \ref{P: cone hom2}. The upper left vertical map is induced by a quasi-isomorphism from a homotopically fine sheaf to an injective sheaf and so is a cohomology isomorphism by \cite[Thm. IV.2.2]{Br}. The map labelled $f$ is a cohomology isomorphism for the same reason, while the two righthand vertical maps are also cohomology isomorphisms by \cite[Thm. IV.2.2]{Br}, being induced by quasi-isomorphisms of injective sheaves. Hence, it follows that the top horizontal map, as desired, is a cohomology isomorphism for $*\leq \bar p(k)$.

This concludes the proof.
\end{proof}

\begin{corollary}
Let $X$ be a paracompact $n$-dimensional  topological stratified pseudomanifold, and let $\mc{G}$ denote a local coefficient system of $R$ modules on $X-\Sigma$. Then $\H^i_c(\mc {P}^*_{\mc G\otimes_R \mc O}) \cong IH^{c}_{n-i}(X;\mc G_0)$. 
\end{corollary}
\begin{proof}
By definition, $\H^i_c(\mc {P}^*)\cong H^i(\Gamma_c(X;\mc I^*)$, where $\mc I^*$ is an injective resolution of $\mc P^*$. Since $\mc P^*$ and $\mc{IS}^*$ are quasi-isomorphic, $\dim_RX<\infty$, and $\mc{IS}^*$ is homotopically fine, \cite[Thm. IV.2.2]{Br} then yields an isomorphism  $\H^i_c(\mc {P}^*)\cong H^i(\Gamma_c(X;\mc {IS}^*))$. So it remains to see that $\Gamma_c(X;\mc {IS}^*)\cong IC_{n-*}^c(X)$.

By \cite[Thm. I.6.2]{Br}, since $IS^*$ has no non-trivial $0$ sections (by Lemma \ref{L: global mono}), $\Gamma_c(X;\mc {IS}^*)$ is isomorphic to the submodule of $IS^*(X)=IC^{\infty}_{n-*}(X)$ consisting of presheaf sections with compact support. But these are exactly the compactly supported (finite) intersection chains $IC_{n-*}^c(X)$. 
\end{proof}

This corollary says that the hypercohomology with compact supports of the Deligne sheaf complex can be computed via finite singular chain intersection homology. In particular, it is a stratum-preserving homotopy invariant:

\begin{corollary}
Let $f:X\to Y$ be a stratum-preserving homotopy equivalence of  paracompact $n$-dimensional stratified topological  pseudomanifolds. Then $\H^i_c(\mc {P}^*_{\mc G\otimes_R \mc O}(Y))\cong \H^i_c(\mc {P}^*_{f^*\mc G\otimes_R f^*\mc O}(X))$. 
\end{corollary}

\subsection{Calculation on PL pseudomanifolds}\label{S: PL}

In this section, we indicate how superperverse intersection homology may be computed on PL pseudomanifolds via direct use of simplicial chains. 

Suppose that $X$ is an $n$-dimensional PL stratified pseudomanifold. We assume $X$ to be second countable, in accordance with the definition of a PL space given in \cite{HUD}. In this case, $X$ is embeddable in some $\R^N$, and, in particular, $X$ is metrizable and thus paracompact \cite[Ch. III]{HUD}. 

Given any fixed triangulation $K$ of $X$ compatible with the stratification (i.e. each skeleton $X^i$ is a subpolyhedron), we can form the simplicial intersection chain complex $\mf{I^{\bar p}_KC}^{\infty}_*(X;\mc G_0)$. In particular, let $\mc G_0$ continue to denote a coefficient system given by $\mc {G}$ on $X-\Sigma=X-X^{n-2}$ and by a $0$ coeficient system on $\Sigma$. Then any simplex $\sigma$ in $K$ can be given a coefficient by lifting the complement in $\sigma$ of $\sigma\cap \Sigma$ to the covering space $\mc G$. Any simplex lying in $\Sigma$ is given a $0$ coefficient. Boundaries of simplices are then computed via the usual formula, taking coefficients given by restricting lifts of faces not contained in $\Sigma$ and discarding (equivalently giving $0$ coefficients to) faces contained in $\Sigma$. An $i$-simplex $\sigma$ is $\bar p$ allowable if $\dim(\sigma\cap X_{n-k})\leq i-k+\bar p(k)$, and an $i$-chain $\xi$ is $\bar p$ allowable if each simplex with non-zero coefficient in $\xi$ and $\bd \xi$ is $\bar p$ allowable. Thus $\mf{I^{\bar p}_KC}^{\infty}_*(X;\mc G_0)$ is a direct generalization of the usual simplicial intersection with closed supports \cite{GM2}. This gives a candidate intersection chain complex, and we claim that its homology gives the same intersection homology modules as obtained via the singular chain complex, at least provided we take a sufficiently fine triangulation.

In particular, as in \cite{GM1} (see also \cite[Ch. I,II]{Bo}), we can form the direct limit of the chain complexes $\mf{I^{\bar p}C}^{\infty}_*(X;\mc G_0)=\lim_K \mf{I^{\bar p}_KC}^{\infty}_*(X;\mc G_0)$ under the direct system given by refinement of triangulations, or, following the appendix by Goresky and MacPherson in \cite{MV86}, we can choose any triangulation compatible with the stratification, subdivide barycentrically  once, and then define intersection chains with respect  to the resulting complex $\mf{I^{\bar p}_{K'}C}^{\infty}_*(X;\mc G_0)$. $\mf{I^{\bar p}C}^{\infty}_*(X;\mc G_0)$ and $ \mf{I^{\bar p}_{K'}C}^{\infty}_*(X;\mc G_0)$ 
will be quasi-isomorphic by the same methods discussed there. Then one can proceed as in \cite{GM2} or \cite{Bo} to produce sheaves of PL intersection chains as the sheafification of the presheaf $U\to \mf{I^{\bar p}C}^{\infty}_*(U;\mc G_0)$. In this case, the sheaf so obtained will be soft, following the same proof as given by Habegger in \cite[\S II.5]{Bo}. The PL computations of the local intersection homology groups over distinguished neighborhoods and deleted distinguished neighborhoods also can proceed as in \cite[\S II]{Bo}, but with sufficient modifications due to the coefficient system $\mc G_0$. We see that we have obtained a sheaf satisfying the  correct axioms for the superperverse intersection homology theory. 

So, by the axiomatic nature of the Deligne sheaf construction and the softness of the sheaf associated to $\mf{I^{\bar p}C}^{\infty}_*(X;\mc G_0)$, we have $\H^*(\mc P^*)\cong H_{n-*}(\mf{I^{\bar p}C}^{\infty}_*(X;\mc G_0))\cong H_{n-*}(\mf{I^{\bar p}_{K'}C}^{\infty}_*(X;\mc G_0))$.

Thus superperverse intersection homology can be computed simplicially on a PL stratified pseudomanifold. Furthermore, if $X$ is compact, then we may use finite chains.

\bibliographystyle{amsplain}
\bibliography{bib}

\end{document}